# MULTIVARIATE QUANTILES AND MULTIPLE-OUTPUT REGRESSION QUANTILES: FROM $L_1$ OPTIMIZATION TO HALFSPACE DEPTH

By Marc Hallin[2,3], Davy Paindaveine[3,4] and Miroslav Šiman[4]

*Université Libre de Bruxelles*

A new multivariate concept of quantile, based on a directional version of Koenker and Bassett's traditional regression quantiles, is introduced for multivariate location and multiple-output regression problems. In their empirical version, those quantiles can be computed efficiently via linear programming techniques. Consistency, Bahadur representation and asymptotic normality results are established. Most importantly, the contours generated by those quantiles are shown to coincide with the classical halfspace depth contours associated with the name of Tukey. This relation does not only allow for efficient depth contour computations by means of parametric linear programming, but also for transferring from the quantile to the depth universe such asymptotic results as Bahadur representations. Finally, linear programming duality opens the way to promising developments in depth-related multivariate rank-based inference.

**1. Introduction: Multivariate quantiles and statistical depth.** In this paper, we propose a definition of multivariate quantiles/multiple-output regression quantiles enjoying all the probabilistic and analytical properties one is generally expecting from a quantile, while exhibiting a very strong and fundamental connection with the concept of halfspace depth. Some of the basic

---

Received December 2008; revised June 2009.

[1]Discussed in 10.1214/09-AOS723A, 10.1214/09-AOS723B and 10.1214/09-AOS723C; rejoinder at 10.1214/09-AOS723REJ.

[2]Supported by the Sonderforschungsbereich "Statistical modelling of nonlinear dynamic processes" (SFB823) of the German Research Foundation (Deutsche Forschungsgemeinschaft) and by a Discovery Grant of the Australian Research Council. Marc Hallin is also member of the Académie Royale de Belgique and of CentER, Tilburg University.

[3]Marc Hallin and Davy Paindaveine are also members of ECORE, the association between CORE and ECARES.

[4]Supported by a Mandat d'Impulsion Scientifique of the Belgian Fonds National de la Recherche Scientifique.

---

ideas of this definition were exposed in an unpublished master thesis by Laine [21], quoted in [16]. In this paper, we carefully revive Laine's ideas, and systematically develop and prove the main properties of the concept he introduced.

A huge literature has been devoted to the problem of extending to a multivariate setting the fundamental one-dimensional concept of quantile; see, for instance, [1, 3–7, 10, 15, 19, 34] and [37] or [33] for a recent survey. An equally huge literature—see [9, 22, 39] and [40] for a comprehensive account—is dealing with the concept of (location) depth. The philosophies underlying those two concepts at first sight are quite different, and even, to some extent, opposite. While quantiles resort to analytical characterizations through inverse distribution functions or $L_1$ optimization, depth often derives from more geometric considerations such as halfspaces, simplices, ellipsoids and projections. Both carry advantages and some drawbacks. Analytical definitions usually bring in efficient algorithms and tractable asymptotics. The geometric ones enjoy attractive equivariance properties and intuitive contents, but their probabilistic study and asymptotics are generally trickier, while their implementation, as a rule, leads to heavy combinatorial algorithms; a highly elegant analytical approach to depth has been proposed in [24], but does not help much in that respect.

Yet, beyond those sharp methodological differences, quantiles and depth obviously exhibit a close conceptional kinship. In the univariate case, all definitions basically agree that the depth of a point $x \in \mathbb{R}$ with respect to a probability distribution P with strictly monotone distribution function $F$ should be $\min(F(x), 1 - F(x))$, so that the only points with depth $d$ are $x_d := F^{-1}(d)$ and $x_{1-d} := F^{-1}(1 - d)$—the quantiles of orders $d$ and $1 - d$, respectively. Starting with dimension two, no such clear and undisputable relation has been established so far—how could there be one, by the way, as long as no clear and undisputable definition of a multivariate quantile has been agreed upon? Bridging the gap between the two concepts thus would allow for transferring to the depth universe the analytical and algorithmic tools of the quantile approach, while sorting out the many candidates for a sound definition of multivariate quantiles. Establishing a relation between the quantile and depth philosophies in $\mathbb{R}^k$, if at all possible, therefore is highly desirable.

An important step in that direction has been made very recently in a paper by Kong and Mizera [20]. Kong and Mizera adopt a very simple and, at first sight, quite natural projection-based definition of quantiles. In that approach, denoting by $\mathbf{u}$ a point on the unit sphere $\mathcal{S}^{k-1}$, a quantile of order $\tau \in (0, 1)$ is either a real number $q_{\mathrm{KM};\tau\mathbf{u}} \in \mathbb{R}$ ($q_{\mathrm{KM};\tau\mathbf{u}}^{(n)}$ in the empirical case), the point $\mathbf{q}_{\mathrm{KM};\tau\mathbf{u}} := q_{\mathrm{KM};\tau\mathbf{u}}\mathbf{u} \in \mathbb{R}^k$ (resp., $\mathbf{q}_{\mathrm{KM};\tau\mathbf{u}}^{(n)}$), or the hyperplane $\pi_{\mathrm{KM};\tau\mathbf{u}}$ (resp., $\pi_{\mathrm{KM};\tau\mathbf{u}}^{(n)}$) orthogonal to $\mathbf{u}$ at $\mathbf{q}_{\mathrm{KM};\tau\mathbf{u}}$ (resp., $\mathbf{q}_{\mathrm{KM};\tau\mathbf{u}}^{(n)}$). The scalar quantity $q_{\mathrm{KM};\tau\mathbf{u}} \in \mathbb{R}$ is defined as the quantile of order $\tau$ of the univariate distribution obtained by projecting P onto the oriented straight line with unit vector $\mathbf{u}$, and therefore derives from purely univariate $L_1$ arguments; see Section 4.3 for details. The resulting quantile contours (the collections, for fixed $\tau$, of $\mathbf{q}_{\mathrm{KM};\tau\mathbf{u}}$'s) do not enjoy the properties (independence with respect to the choice of an origin, affine-equivariance, nestedness, etc.) one is expecting from a quantile concept. However, somewhat surprisingly, the envelopes of these contours—namely, the inner regions characterized by the (infinite) fixed-$\tau$ collections of $\pi_{\mathrm{KM};\tau\mathbf{u}}$'s (resp., $\pi_{\mathrm{KM};\tau\mathbf{u}}^{(n)}$'s)—coincide with Tukey's halfspace depth regions, which provides a most interesting, though somewhat indirect, conceptual bridge between the two concepts.

Our quantiles also are associated with unit vectors $\mathbf{u} \in \mathcal{S}^{k-1}$, hence also are directional quantiles. However, instead of projecting onto the straight line defined by $\mathbf{u}$, we stay in a $k$-dimensional setting, where $\mathbf{u}$ simply indicates the reference "vertical" direction for a regression quantile construction in the Koenker and Bassett [18] style. As in [18], our quantiles thus are hyperplanes $\pi_{\tau\mathbf{u}}$ ($\pi_{\tau\mathbf{u}}^{(n)}$ in the empirical case); in contrast with $\pi_{\mathrm{KM};\tau\mathbf{u}}$ and $\pi_{\mathrm{KM};\tau\mathbf{u}}^{(n)}$, however, the fixed-$\mathbf{u}$ collections of $\pi_{\tau\mathbf{u}}$'s and $\pi_{\tau\mathbf{u}}^{(n)}$'s are not collections of parallel hyperplanes all orthogonal to $\mathbf{u}$. Whereas projection quantiles only involve univariate $L_1$ arguments, ours indeed rely on fully $k$-dimensional $L_1$ optimization. As shown in Section 4, the inner regions characterized by the fixed-$\tau$ collections of $\pi_{\tau\mathbf{u}}$'s (resp., $\pi_{\tau\mathbf{u}}^{(n)}$'s) also coincide with Tukey's halfspace depth regions. Contrary to Kong and Mizera's, however, the $\pi_{\tau\mathbf{u}}$ quantile hyperplanes do enjoy all the desirable properties of a well-behaved quantile concept. And, in the empirical case, our quantile hyperplanes and the faces of Tukey's (polyhedral) depth contours essentially coincide, in the sense that the latter constitute a (finite) subcollection of the finite collection of $\pi_{\tau\mathbf{u}}^{(n)}$'s, itself a finite subcollection of the *infinite* collection of $\pi_{\mathrm{KM};\tau\mathbf{u}}^{(n)}$'s.

From their $L_1$ definitions, the $\pi_{\tau\mathbf{u}}^{(n)}$'s and $\pi_{\mathrm{KM};\tau\mathbf{u}}^{(n)}$'s both inherit a probabilistic interpretation allowing for tractable asymptotics: consistency, Bahadur representations and asymptotic normality. From their relation to depth, the resulting contours acquire a series of nice geometric properties such as convexity, nestedness and affine-equivariance; and, since empirical Tukey depth contours fully characterize the empirical distribution (see [35]), our quantile contours (as well as Kong and Mizera's) also do. Above all, our quantiles receive the important benefits of linear programming algorithms, which thereby automatically transfer to depth, hence—indirectly, though (see [26])—also to the Kong and Mizera concept. Moreover, both concepts readily generalize to the regression setting, yielding nested polyhedral regions wrapping, up to the classical quantile crossings, a median or deepest regression hypertube (see [26] for a detailed comparison of our regression quantile hypertubes and those resulting from the Kong and Mizera approach). This extends to the multiple-output context the celebrated single-output Koenker and Bassett concept of regression quantiles. Conversely, as indicated in [26] it also leads to a concept of *multiple-output regression halfspace depth*; that depth concept, however, has the nature of a *point regression depth*, hence is distinct from the Rousseeuw and Hubert regression depth concept (see [29]), which is a *hyperplane depth* concept. A constrained optimization form of the definition of $\pi_{\tau\mathbf{u}}^{(n)}$ also allows for computing Lagrange multipliers with most interesting statistical applications. Finally, by resorting to classical linear programming duality, a concept of directional regression rank scores, allowing for multivariate versions of the methods developed in [12], naturally comes into the picture.

From an applied perspective, the possibility of computing Tukey depth contours via parametric linear programming is not a small advantage. The complexity of computing the depth of a given point is $O(n^{k-1}\log n)$, with algorithms by Rousseeuw and Ruts [28] for $k = 2$ and Rousseeuw and Struyf [31] for general $k$. The best known algorithm for computing all depth contours has complexity $O(n^2)$ (see [23]) in dimension $k = 2$. To the best of our knowledge, no exact implementable algorithm is available so far for $k > 2$. Our approach allows for higher values of $k$, and we could easily run our algorithms in dimension $k = 5$, for a few hundred observations.

The paper is organized as follows. Section 2 introduces the definitions and main notation to be used throughout. In Section 3, we study the main properties of the new quantiles: from their directional quantile nature, they inherit subgradient characterizations (Section 3.1), equivariance properties (Section 3.2), and quantile-like asymptotics—strong consistency, Bahadur representation and asymptotic normality (Section 3.3). In Section 4, we establish the equivalence of the quantile contours, thus obtained with the more traditional halfspace (or Tukey) depth contours, as well as their relation to the recent results by Kong and Mizera [20] and Wei [37]. Section 5 is devoted to the computational aspects of our multivariate quantiles, and Section 6 to their extension to a multiple-output regression context. A brief application to real data is discussed in Section 7. Section 8 concludes with some perspectives for future research. Proofs are collected in the Appendix.

**2. Definition and notation.** Consider the $k$-variate random vector $\mathbf{Z} := (Z_1, \ldots, Z_k)'$. The multivariate quantiles we are proposing are *directional* objects—more precisely, $(k-1)$-dimensional hyperplanes indexed by vectors $\boldsymbol{\tau}$ ranging over the open unit ball (deprived of the origin) $\mathcal{B}^k := \{\mathbf{z} \in \mathbb{R}^k : 0 < \|\mathbf{z}\| < 1\}$ of $\mathbb{R}^k$. This directional index $\boldsymbol{\tau}$ naturally factorizes into $\boldsymbol{\tau} =: \tau\mathbf{u}$, where $\tau = \|\boldsymbol{\tau}\| \in (0, 1)$ and $\mathbf{u} \in \mathcal{S}^{k-1} := \{\mathbf{z} \in \mathbb{R}^k : \|\mathbf{z}\| = 1\}$. Denoting by $\boldsymbol{\Gamma}_{\mathbf{u}}$ an arbitrary $k \times (k-1)$ matrix of unit vectors such that $(\mathbf{u} \vdots \boldsymbol{\Gamma}_{\mathbf{u}})$ constitutes an orthonormal basis of $\mathbb{R}^k$, we define the $\boldsymbol{\tau}$-quantile of $\mathbf{Z}$ as the regression $\tau$-quantile obtained (in the traditional Koenker and Bassett [18] sense) when regressing $\mathbf{Z_u} := \mathbf{u}'\mathbf{Z}$ on the marginals of $\mathbf{Z_u^\perp} := \boldsymbol{\Gamma}_\mathbf{u}'\mathbf{Z}$ and a constant term: the vector $\mathbf{u}$ therefore indicates the direction of the "vertical" axis in the regression, while $\boldsymbol{\Gamma}_\mathbf{u}$ simply provides an orthonormal basis of the vector space orthogonal to $\mathbf{u}$. More precisely, denoting by $x \mapsto \rho_\tau(x) := x(\tau - \mathbb{I}_{[x<0]})$ the usual $\tau$-quantile *check function*, we adopt the following definition.

DEFINITION 2.1. The $\boldsymbol{\tau}$-quantile of $\mathbf{Z}$ ($\boldsymbol{\tau} =: \tau\mathbf{u} \in \mathcal{B}^k$) is any element of the collection $\Pi_{\boldsymbol{\tau}}$ of hyperplanes $\pi_{\boldsymbol{\tau}} := \{\mathbf{z} \in \mathbb{R}^k : \mathbf{u}'\mathbf{z} = \mathbf{b}_{\boldsymbol{\tau}}'\boldsymbol{\Gamma}_\mathbf{u}'\mathbf{z} + a_{\boldsymbol{\tau}}\}$ such that

$$(2.1) \quad \begin{aligned} (a_{\boldsymbol{\tau}}, \mathbf{b}_{\boldsymbol{\tau}}')' &\in \underset{(a,\mathbf{b}')' \in \mathbb{R}^k}{\arg\min} \ \Psi_{\boldsymbol{\tau}}(a, \mathbf{b}) \\ &\text{where } \Psi_{\boldsymbol{\tau}}(a, \mathbf{b}) := \mathrm{E}[\rho_\tau(\mathbf{Z_u} - \mathbf{b}'\mathbf{Z_u^\perp} - a)]. \end{aligned}$$

This definition tacitly requires the existence, for $\mathbf{Z}$, of finite first-order moments: see the comment below. For the sake of notational simplicity, quantiles, here and in the sequel, are associated with a random vector $\mathbf{Z}$, though they actually are attributes of $\mathbf{Z}$'s probability distribution P.

Definition 2.1 clearly extends the traditional univariate one. For $k = 1$, indeed, hyperplanes of dimension $k - 1$ are simply points, $\mathcal{B}^k$ reduces to $(-1, 0) \cup (0, 1)$ and $\pi_{\boldsymbol{\tau}}$ to a "classical" quantile, of order $1 - \|\boldsymbol{\tau}\|$ ($\boldsymbol{\tau}$ pointing to the left) or $\|\boldsymbol{\tau}\|$ ($\boldsymbol{\tau}$ pointing to the right). This couple of quantiles constitutes (for $k = 1$) a *quantile contour*, indicating that a sensible relation between depth and quantiles should associate depth contours with contour-valued rather than with point-valued quantiles.

Note that the quantile hyperplanes $\pi_{\boldsymbol{\tau}}$ and the "intercepts" $a_{\boldsymbol{\tau}}$ are well defined in the sense that they only depend on $\boldsymbol{\tau}$, not on the coordinate system associated with the (arbitrary) choice of $\boldsymbol{\Gamma}_\mathbf{u}$. However, the "slope" coefficients $\mathbf{b}_{\boldsymbol{\tau}} = \mathbf{b}_{\boldsymbol{\tau}}(\boldsymbol{\Gamma}_\mathbf{u})$ do depend on $\boldsymbol{\Gamma}_\mathbf{u}$, a dependence we do not stress in the notation unless really necessary.

Each quantile hyperplane $\pi_{\boldsymbol{\tau}}$ [each element $(a_{\boldsymbol{\tau}}, \mathbf{b}_{\boldsymbol{\tau}}')'$ of $\arg\min_{(a,\mathbf{b}')' \in \mathbb{R}^k} \Psi_{\boldsymbol{\tau}}(a, \mathbf{b})$] characterizes a lower (open) *quantile halfspace*

$$(2.2) \quad H_{\boldsymbol{\tau}}^- = H_{\boldsymbol{\tau}}^-(a_{\boldsymbol{\tau}}, \mathbf{b}_{\boldsymbol{\tau}}) := \{\mathbf{z} \in \mathbb{R}^k : \mathbf{u}'\mathbf{z} < \mathbf{b}_{\boldsymbol{\tau}}'\boldsymbol{\Gamma}_\mathbf{u}'\mathbf{z} + a_{\boldsymbol{\tau}}\}$$

and an upper (closed) *quantile halfspace*

$$(2.3) \quad H_{\boldsymbol{\tau}}^+ = H_{\boldsymbol{\tau}}^+(a_{\boldsymbol{\tau}}, \mathbf{b}_{\boldsymbol{\tau}}) := \{\mathbf{z} \in \mathbb{R}^k : \mathbf{u}'\mathbf{z} \geq \mathbf{b}_{\boldsymbol{\tau}}'\boldsymbol{\Gamma}_\mathbf{u}'\mathbf{z} + a_{\boldsymbol{\tau}}\}.$$

As already mentioned, Definition 2.1 requires $\mathbf{Z}$ to have finite first-order moments. Actually, modifying (2.1) into $(a_{\boldsymbol{\tau}}, \mathbf{b}_{\boldsymbol{\tau}}')' \in \arg\min_{(a,\mathbf{b}')' \in \mathbb{R}^k}(\Psi_{\boldsymbol{\tau}}(a, \mathbf{b}) - \Psi_{\boldsymbol{\tau}}(0, \mathbf{0}))$ has no impact on $\pi_{\boldsymbol{\tau}}$, while allowing to relax the moment condition on $\mathbf{Z_u}$; finite first-order moments, however, still are required for $\mathbf{Z_u^\perp}$. When $\mathbf{u}$ ranges over $\mathcal{S}^{k-1}$—for instance, when defining quantile contours—we need finite first-order moments for all $\mathbf{Z_u^\perp}$'s, hence for $\mathbf{Z}$ itself. For the sake of simplicity, we often adopt the following assumption in the sequel.

ASSUMPTION (A).   The distribution of the random vector $\mathbf{Z}$ is absolutely continuous with respect to the Lebesgue measure on $\mathbb{R}^k$, with a density ($f$, say) that has connected support, and admits finite first-order moments.

The minimization problem (2.1) may have multiple solutions, yielding distinct hyperplanes $\pi_{\boldsymbol{\tau}}$. This, however, does not occur under Assumption (A), as shown in the following result, which is a particular case of Theorem 2.1 in [26].

PROPOSITION 2.1.   *Let Assumption (A) hold. Then, for any $\boldsymbol{\tau} \in \mathcal{B}^k$, the minimizer $(a_{\boldsymbol{\tau}}, \mathbf{b}'_{\boldsymbol{\tau}})'$ in (2.1), hence also the resulting quantile hyperplane $\pi_{\boldsymbol{\tau}}$, is unique.*

The family of hyperplanes $\Pi = \{\pi_{\boldsymbol{\tau}} : \boldsymbol{\tau} = \tau\mathbf{u} \in \mathcal{B}^k\}$ can be considered from two different points of view. The *directional* point of view, associated with the fixed-$\mathbf{u}$ subfamilies $\Pi_{\mathbf{u}} := \{\pi_{\boldsymbol{\tau}} : \boldsymbol{\tau} = \tau\mathbf{u}, \tau \in (0, 1)\}$ is the one emphasized so far in the definition, and provides, for each $\mathbf{u}$, the usual interpretation of a collection of regression quantile hyperplanes. Another point of view is associated with the fixed-$\tau$ subfamilies $\Pi_\tau := \{\pi_{\boldsymbol{\tau}} : \boldsymbol{\tau} = \tau\mathbf{u}, \mathbf{u} \in \mathcal{S}^{k-1}\}$, which generate *quantile contours*: this point of view is developed in Section 4.

Before turning to the empirical version of our quantiles, let us present an alternative (but strictly equivalent) definition of our $\boldsymbol{\tau}$-quantiles, based on a *constrained* optimization formulation.

DEFINITION 2.2.   The $\boldsymbol{\tau}$-quantile of $\mathbf{Z}$ ($\boldsymbol{\tau} =: \tau\mathbf{u} \in \mathcal{B}^k$) is any element of the collection $\Pi_{\boldsymbol{\tau}}$ of hyperplanes $\pi_{\boldsymbol{\tau}} := \{\mathbf{z} \in \mathbb{R}^k : \mathbf{c}'_{\boldsymbol{\tau}}\mathbf{z} = a_{\boldsymbol{\tau}}\}$ such that

$$(2.4) \qquad (a_{\boldsymbol{\tau}}, \mathbf{c}'_{\boldsymbol{\tau}})' \in \underset{(a, \mathbf{c}')' \in \mathcal{M}_{\mathbf{u}}}{\arg\min} \ \Psi^c_\tau(a, \mathbf{c}),$$

where $\Psi^c_\tau(a, \mathbf{c}) := \mathrm{E}[\rho_\tau(\mathbf{c}'\mathbf{Z} - a)]$ and $\mathcal{M}_{\mathbf{u}} := \{(a, \mathbf{c}')' \in \mathbb{R}^{k+1} : \mathbf{u}'\mathbf{c} = 1\}$.

Clearly, if $(a_{\boldsymbol{\tau}}, \mathbf{b}'_{\boldsymbol{\tau}})'$ is a minimizer of (2.1), then $(a_{\boldsymbol{\tau}}, \mathbf{c}'_{\boldsymbol{\tau}})' := (a_{\boldsymbol{\tau}}, (\mathbf{u} - \boldsymbol{\Gamma_u}\mathbf{b}_{\boldsymbol{\tau}})')'$ minimizes the objective function in (2.4); conversely, for any minimizer $(a_{\boldsymbol{\tau}}, \mathbf{c}'_{\boldsymbol{\tau}})'$ of (2.4), $(a_{\boldsymbol{\tau}}, \mathbf{b}'_{\boldsymbol{\tau}})' := (a_{\boldsymbol{\tau}}, (-\boldsymbol{\Gamma'_u}\mathbf{c}_{\boldsymbol{\tau}})')'$ minimizes the objective function in (2.1). The two definitions thus coincide; in particular, the lower and upper quantile halfspaces $\{\mathbf{z} \in \mathbb{R}^k : \mathbf{c}'_{\boldsymbol{\tau}}\mathbf{z} < a_{\boldsymbol{\tau}}\}$ and $\{\mathbf{z} \in \mathbb{R}^k : \mathbf{c}'_{\boldsymbol{\tau}}\mathbf{z} \geq a_{\boldsymbol{\tau}}\}$ associated with the quantile hyperplanes of Definition 2.2 coincide with those in (2.2) and (2.3), and therefore, depending on the context, the notation $H_{\tau}^{\pm}(a_{\tau}, \mathbf{b}_{\tau})$, $H_{\tau}^{\pm}(a_{\tau}, \mathbf{c}_{\tau})$, or simply $H_{\tau}^{\pm}$ will be used indifferently. Definitions 2.1 and 2.2 both have advantages and, in the sequel, we use them both. Definition 2.1 is preferred in this section since it carries all the intuitive contents of our concept; the advantages of Definition 2.2, of an analytical nature, will appear more clearly in Sections 3.1 and 5.

The empirical versions of our quantile hyperplanes and the corresponding lower and upper quantile halfspaces naturally follow as sample analogs of the population concepts. To be more specific, let $\mathbf{Z}^{(n)} := (\mathbf{Z}_1, \ldots, \mathbf{Z}_n)$ be an $n$-tuple $(n > k)$ of $k$-dimensional random vectors: we define the *empirical $\boldsymbol{\tau}$-quantile* of $\mathbf{Z}^{(n)}$ as any element of the collection $\Pi_{\boldsymbol{\tau}}^{(n)}$ of hyperplanes $\pi_{\boldsymbol{\tau}}^{(n)} := \{\mathbf{z} \in \mathbb{R}^k : \mathbf{u}'\mathbf{z} = \mathbf{b}_{\boldsymbol{\tau}}^{(n)'} \boldsymbol{\Gamma}_{\mathbf{u}}' \mathbf{z} + a_{\boldsymbol{\tau}}^{(n)}\}$ such that (with obvious notation)

$$(2.5) \quad (a_{\boldsymbol{\tau}}^{(n)}, \mathbf{b}_{\boldsymbol{\tau}}^{(n)'})' \in \underset{(a, \mathbf{b}')' \in \mathbb{R}^k}{\arg\min} \ \Psi_{\boldsymbol{\tau}}^{(n)}(a, \mathbf{b})$$

$$\text{with } \Psi_{\boldsymbol{\tau}}^{(n)}(a, \mathbf{b}) := \frac{1}{n} \sum_{i=1}^{n} \rho_{\tau}(\mathbf{Z}_{i\mathbf{u}} - \mathbf{b}'\mathbf{Z}_{i\mathbf{u}}^{\perp} - a),$$

or equivalently, of hyperplanes $\pi_{\boldsymbol{\tau}}^{(n)} := \{\mathbf{z} \in \mathbb{R}^k : \mathbf{c}_{\boldsymbol{\tau}}^{(n)'} \mathbf{z} = a_{\boldsymbol{\tau}}^{(n)}\}$ such that

$$(2.6) \quad (a_{\boldsymbol{\tau}}^{(n)}, \mathbf{c}_{\boldsymbol{\tau}}^{(n)'})' \in \underset{(a, \mathbf{c}')' \in \mathcal{M}_{\mathbf{u}}}{\arg\min} \ \Psi_{\boldsymbol{\tau}}^{c(n)}(a, \mathbf{c})$$

$$\text{with } \Psi_{\boldsymbol{\tau}}^{c(n)}(a, \mathbf{c}) := \frac{1}{n} \sum_{i=1}^{n} \rho_{\tau}(\mathbf{c}'\mathbf{Z}_i - a)$$

(no moment assumption is required here). These empirical quantiles—which for given $\mathbf{u}$ clearly coincide with the Koenker and Bassett [18] hyperplanes in the coordinate system $(\mathbf{u} \vdots \boldsymbol{\Gamma}_{\mathbf{u}})$—allow for defining, in an obvious way, the empirical analogs $H_{\tau}^{(n)-}$ and $H_{\tau}^{(n)+}$ of the lower and upper quantile halfspaces in (2.2) and (2.3); see Figures 1 and 2 for an illustration.

Of course, empirical distributions are inherently discrete, and empirical $\boldsymbol{\tau}$-quantiles and halfspaces in general are not uniquely defined. However, the minimizers of (2.5) [equivalently, of (2.6)], for given $\boldsymbol{\tau}$, are "close to each other," in the sense that the set of minimizers is convex—hence, connected (this readily follows from the fact that the objective functions are convex); this set is shrinking, as $n \to \infty$, to a single point which corresponds to the uniquely defined population quantile, provided that the following assumption is fulfilled (see the asymptotic results of Section 3.3 for details).

ASSUMPTION ($A_n$).   The observations $\mathbf{Z}_i$, $i = 1, \ldots, n$ are i.i.d. with a common distribution satisfying Assumption (A).

Finally, note that, since the empirical versions of our quantiles, for given $\mathbf{u}$, are defined as standard single-output quantile regression hyperplanes, they inherit the linear programming features of the Koenker–Bassett theory. This certainly is one of the most important and attractive properties of the proposed quantiles; see Section 5 for details.

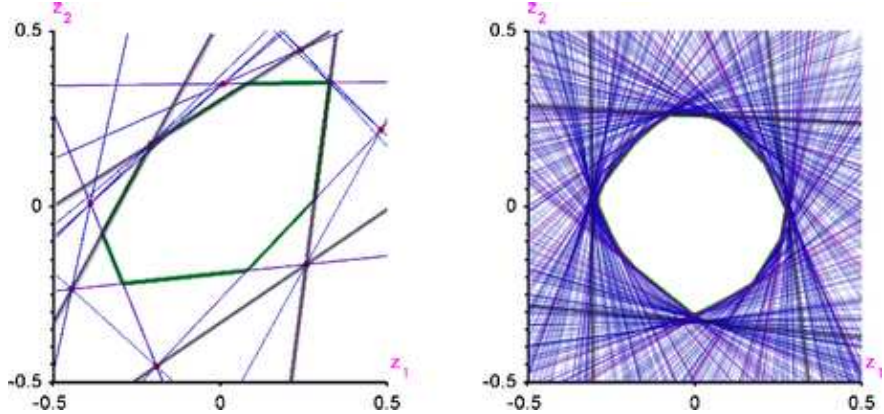

Fig. 1. *The left plot contains $n = 9$ (red) points drawn from $U([-0.5, 0.5]^2)$, the centered bivariate uniform distribution over the unit square, and provides all $\tau$-quantile hyperplanes for $\tau = 0.2$. These hyperplanes define a polygonal central region (green contour) which, in Section 4, is shown to coincide with a Tukey depth region. The quantile hyperplanes contributing an edge to the polygonal central region are shown in magenta; those associated with the four semiaxial directions in black; all other ones in blue. The right plot provides the same information for $n = 499$ (invisible) points drawn from the same population distribution.*

**3. Multivariate quantiles as directional quantiles.** In this section, we describe the "directional" properties of our quantiles. We first derive and dis-cuss the subgradient conditions associated with the optimization problems in Definitions 2.1 and 2.2, then state the strong equivariance properties of our empirical quantiles, and finally present some asymptotic results.

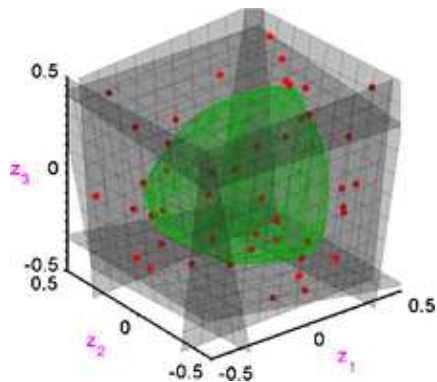

Fig. 2. *This plot provides six $\tau$-quantile hyperplanes (in black) in the semiaxial directions for $\tau = 0.1$, computed from $n = 49$ (red) points drawn from $U([-0.5, 0.5]^3)$, the centered trivariate uniform distribution over the unit cube, along with the corresponding central (Tukey depth) region (in green).*

3.1. *Subgradient conditions.* Under Assumption (A), the objective function $\Psi_{\boldsymbol{\tau}}$ appearing in Definition 2.1 is convex and continuously differentiable on $\mathbb{R}^k$. Therefore, our population $\boldsymbol{\tau}$-quantiles can be equivalently defined as the collection of hyperplanes associated with the solutions $(a_{\boldsymbol{\tau}}, \mathbf{b}'_{\boldsymbol{\tau}})'$ of the system of equations

$$\text{(3.1)} \qquad \operatorname{grad}_{(a, \mathbf{b}')'} \Psi_{\boldsymbol{\tau}}(a, \mathbf{b}) = \mathbf{0}$$

(see Sections 2.2.1 and 2.2.2 in [16]). These hyperplanes thus are characterized by the relations

$$
\begin{aligned}
\text{(3.2a)} \quad 0 &= (\partial_a \Psi_{\boldsymbol{\tau}}(a, \mathbf{b}))_{(a_{\boldsymbol{\tau}}, \mathbf{b}'_{\boldsymbol{\tau}})'} = \mathrm{P}[\mathbf{u}'\mathbf{Z} < \mathbf{b}'_{\boldsymbol{\tau}}\boldsymbol{\Gamma}'_{\mathbf{u}}\mathbf{Z} + a_{\boldsymbol{\tau}}] - \tau \\
&= \mathrm{P}[\mathbf{Z} \in H^-_{\boldsymbol{\tau}}(a_{\boldsymbol{\tau}}, \mathbf{b}_{\boldsymbol{\tau}})] - \tau,
\end{aligned}
$$

$$
\text{(3.2b)} \quad \mathbf{0} = (\operatorname{grad}_{\mathbf{b}} \Psi_{\boldsymbol{\tau}}(a, \mathbf{b}))_{(a_{\boldsymbol{\tau}}, \mathbf{b}'_{\boldsymbol{\tau}})'} = -\tau \mathrm{E}[\boldsymbol{\Gamma}'_{\mathbf{u}}\mathbf{Z}] + \mathrm{E}[\boldsymbol{\Gamma}'_{\mathbf{u}}\mathbf{Z}\mathbb{I}_{[\mathbf{Z} \in H^-_{\boldsymbol{\tau}}(a_{\boldsymbol{\tau}}, \mathbf{b}_{\boldsymbol{\tau}})]}].
$$

Clearly, relation (3.2a) provides our multivariate $\boldsymbol{\tau}$-quantiles with a natural probabilistic interpretation, as it keeps the probability of their lower halfspaces equal to $\tau(=\|\boldsymbol{\tau}\|)$. As for relation (3.2b), it can be rewritten as

$$\text{(3.3)} \qquad \boldsymbol{\Gamma}'_{\mathbf{u}}\left[\frac{1}{1-\tau}\mathrm{E}[\mathbf{Z}\mathbb{I}_{[\mathbf{Z} \in H^+_{\boldsymbol{\tau}}]}] - \frac{1}{\tau}\mathrm{E}[\mathbf{Z}\mathbb{I}_{[\mathbf{Z} \in H^-_{\boldsymbol{\tau}}]}]\right] = \mathbf{0},$$

which—combined with (3.2a)—shows that the straight line through the probability mass centers $\frac{1}{\tau}\mathrm{E}[\mathbf{Z}\mathbb{I}_{[\mathbf{Z} \in H^-_{\boldsymbol{\tau}}]}]$ and $\frac{1}{1-\tau}\mathrm{E}[\mathbf{Z}\mathbb{I}_{[\mathbf{Z} \in H^+_{\boldsymbol{\tau}}]}]$ of the lower and upper $\boldsymbol{\tau}$-quantile halfspaces is parallel to $\mathbf{u}(:= \boldsymbol{\tau}/\tau)$. Note moreover that, quite trivially,

$$(1-\tau)\left(\frac{1}{1-\tau}\mathrm{E}[\mathbf{Z}\mathbb{I}_{[\mathbf{Z} \in H^+_{\boldsymbol{\tau}}]}]\right) + \tau\left(\frac{1}{\tau}\mathrm{E}[\mathbf{Z}\mathbb{I}_{[\mathbf{Z} \in H^-_{\boldsymbol{\tau}}]}]\right) = \mathrm{E}[\mathbf{Z}],$$

so that the overall probability mass center also belongs to the same straight line.

Now consider the gradient conditions associated with Definition 2.2, which state that $(a_{\boldsymbol{\tau}}, \mathbf{c}'_{\boldsymbol{\tau}}, \lambda_{\boldsymbol{\tau}})'$ are solutions of the system

$$
\begin{aligned}
\text{(3.4)} \qquad &\operatorname{grad}_{(a, \mathbf{c}', \lambda)'} L_{\boldsymbol{\tau}}(a, \mathbf{c}, \lambda) = \mathbf{0} \\
&\text{with } L_{\boldsymbol{\tau}}(a, \mathbf{c}, \lambda) := \Psi^c_{\boldsymbol{\tau}}(a, \mathbf{c}) - \lambda(\mathbf{u}'\mathbf{c} - 1)
\end{aligned}
$$

(the Lagrangian function of the problem). Equivalently [indeed, the only points in $\mathbb{R}^{k+2}$ where $(a, \mathbf{c}', \lambda)' \mapsto L_{\boldsymbol{\tau}}(a, \mathbf{c}, \lambda)$ is not continuously differentiable are of the form $(0, \mathbf{0}', \lambda)'$, hence cannot be associated with a minimum of (2.4)], the latter gradient conditions can be rewritten as

$$
\text{(3.5a)} \quad 0 = (\partial_a L_\tau(a, \mathbf{c}, \lambda))_{(a_\tau, \mathbf{c}'_\tau, \lambda_\tau)'}
$$
$$
= \mathrm{P}[\mathbf{c}'_\tau \mathbf{Z} < a_\tau] - \tau = \mathrm{P}[\mathbf{Z} \in H^-_\tau(a_\tau, \mathbf{c}_\tau)] - \tau,
$$

$$
\text{(3.5b)} \quad \mathbf{0} = (\mathrm{grad}_{\mathbf{c}} L_\tau(a, \mathbf{c}, \lambda))_{(a_\tau, \mathbf{c}'_\tau, \lambda_\tau)'} = \tau \mathrm{E}[\mathbf{Z}] - \mathrm{E}[\mathbf{Z}\mathbb{I}_{[\mathbf{Z} \in H^-_\tau(a_\tau, \mathbf{c}_\tau)]}] - \lambda_\tau \mathbf{u},
$$

$$
\text{(3.5c)} \quad 0 = (\partial_\lambda L_\tau(a, \mathbf{c}, \lambda))_{(a_\tau, \mathbf{c}'_\tau, \lambda_\tau)'} = 1 - \mathbf{u}' \mathbf{c}_\tau.
$$

For such a constrained optimization problem, gradient conditions in general are necessary but not sufficient. In this case, however, note that premultiplying both sides of (3.5b) by $\mathbf{\Gamma}'_{\mathbf{u}}$ yields (3.2b), which clearly implies that, disregarding the Lagrange multiplier $\lambda_\tau$ and (3.5c) to focus on (the coefficients of) the quantile hyperplane $\pi_\tau$, the necessary conditions (3.5a) and (3.5b) are no weaker than the necessary and sufficient ones in (3.2a) and (3.2b), hence are necessary and sufficient, too.

The gradient conditions (3.4) associated with Definition 2.2 are, in a sense, richer than those (3.1) associated with the original definition of our quantiles, which is actually one of the main reasons why we also consider that alternative definition. Indeed, (3.5b), which can be rewritten as

$$
\text{(3.6)} \qquad \frac{1}{1-\tau} \mathrm{E}[\mathbf{Z}\mathbb{I}_{[\mathbf{Z} \in H^+_\tau]}] - \frac{1}{\tau} \mathrm{E}[\mathbf{Z}\mathbb{I}_{[\mathbf{Z} \in H^-_\tau]}] = \frac{\lambda_\tau}{\tau(1-\tau)} \mathbf{u},
$$

is more informative than (3.2b)–(3.3), and clarifies the role of the Lagrange multiplier $\lambda_\tau$. Such a multiplier, which in general only measures the impact of the boundary constraint [in this case, the constraint (3.5c)], here appears as a functional that is potentially useful for testing (central, elliptical, or spherical) symmetry or for measuring directional outlyingness and tail behavior of the distribution; see Section 8. Moreover, premultiplying (3.5b) with $\mathbf{c}'_\tau$ yields $\lambda_\tau(\mathbf{c}'_\tau \mathbf{u}) = \mathrm{E}[(\tau - \mathbb{I}_{[\mathbf{c}'_\tau \mathbf{Z} - a_\tau < 0]})\mathbf{c}'_\tau \mathbf{Z}]$, that is, by using (3.5a) and (3.5c),

$$
\text{(3.7)} \qquad\qquad \lambda_\tau = \Psi^c_\tau(a_\tau, \mathbf{c}_\tau),
$$

so that $\lambda_\tau$ is nothing but the minimum achieved in (2.4) [equivalently, in (2.1)].

The sample objective functions $\Psi^{(n)}_\tau(a, \mathbf{b})$ and $\Psi^{c(n)}_\tau(a, \mathbf{c})$ in (2.5) and (2.6) are not continuously differentiable. They however have directional derivatives in all directions, which can be used to formulate fixed-$\mathbf{u}$ subgradient conditions for the empirical $\boldsymbol{\tau}$-quantiles, $\boldsymbol{\tau} = \tau\mathbf{u}$. Focusing first on the constrained optimization problem (2.6), it is easy to show that the coefficients $(a^{(n)}_\tau, \mathbf{c}^{(n)\prime}_\tau)'$ and the corresponding Lagrange multiplier $\lambda^{(n)}_\tau$ of any empirical $\boldsymbol{\tau}$-quantile $\pi^{(n)}_\tau = \{\mathbf{z} \in \mathbb{R}^k : \mathbf{c}^{(n)\prime}_\tau \mathbf{z} = a^{(n)}_\tau\}$ must satisfy (letting

$r_{i\tau}^{(n)} := \mathbf{c}_\tau^{(n)\prime}\mathbf{Z}_i - a_\tau^{(n)}, \; i = 1, \ldots, n)$

(3.8a) $$\frac{1}{n}\sum_{i=1}^{n}\mathbb{I}_{[r_{i\tau}^{(n)}<0]} \leq \tau \leq \frac{1}{n}\sum_{i=1}^{n}\mathbb{I}_{[r_{i\tau}^{(n)}\leq0]},$$

(3.8b) $$-\frac{1}{n}\sum_{i=1}^{n}\mathbf{Z}_i^-\mathbb{I}_{[r_{i\tau}^{(n)}=0]} \leq \tau\left[\frac{1}{n}\sum_{i=1}^{n}\mathbf{Z}_i\right] - \left[\frac{1}{n}\sum_{i=1}^{n}\mathbf{Z}_i\mathbb{I}_{[r_{i\tau}^{(n)}<0]}\right] - \lambda_\tau^{(n)}\mathbf{u}$$
$$\leq \frac{1}{n}\sum_{i=1}^{n}\mathbf{Z}_i^+\mathbb{I}_{[r_{i\tau}^{(n)}=0]} \quad\text{and}$$

(3.8c) $$0 = 1 - \mathbf{u}'\mathbf{c}_\tau^{(n)},$$

where $\mathbf{z}^+ := (\max(z_1,0),\ldots,\max(z_k,0))'$ and $\mathbf{z}^- := (-\min(z_1,0),\ldots, -\min(z_k,0))'$. These necessary conditions are obtained by imposing that, at $(a_\tau^{(n)},\mathbf{c}_\tau^{(n)\prime},\lambda_\tau^{(n)})'$, directional derivatives in each of the $2(k+2)$ semi-axial directions of the $(a,\mathbf{c}',\lambda)'$-space be nonnegative for $(a,\mathbf{c}')'$ and zero for $\lambda$.

For $n \gg k$, we clearly may interpret (3.8a) and (3.8c) as an approximate version of their population analogs (3.5a) and (3.5b), roughly with the same consequences [the condition (3.8c) simply restates our boundary constraint]. More specifically, (3.8a) indicates that

(3.9) $$\frac{N}{n} \leq \tau \leq \frac{N+Z}{n}, \quad\text{hence}\quad \frac{P}{n} \leq 1-\tau \leq \frac{P+Z}{n},$$

where $N$, $P$ and $Z$ are the numbers of negative, positive and zero values, respectively, in the residual series $r_{i\tau}^{(n)}$, $i = 1, \ldots, n$. This implies that, for noninteger values of $n\tau$, empirical $\boldsymbol{\tau}$-quantile hyperplanes have to go through some of the $\mathbf{Z}_i$'s. Actually, if the data points are in general position [which of course holds with probability one under Assumption (A$_n$)], there exists a sample $\boldsymbol{\tau}$-quantile hyperplane $\pi_\tau^{(n)}$ which fits exactly $k$ observations; (3.9) then holds with $Z = k$ (see Sections 2.2.1 and 2.2.2 of [16]). Note that the inequalities in (3.8a) and (3.8c) [hence, also in (3.9)] must be strict if the sample $\boldsymbol{\tau}$-quantile is to be uniquely defined. Finally, as we will see in (5.2) below, the value of $\lambda_\tau^{(n)}$, parallel to the population case, is the minimal one that can be achieved in (2.6), hence also in (2.5).

For the unconstrained definition of our empirical quantiles in (2.5), necessary and sufficient subgradient conditions can be obtained by applying Theorem 2.1 of [16], since (2.5) is nothing but a standard single-output quantile regression optimization problem. Assuming that the data points are in general position and defining, for any $k$-tuple of indices $h = (i_1, \ldots, i_k)$, $1 \leq i_1 < \cdots < i_k \leq n$,

(3.10) $$\mathbb{Y}_\mathbf{u}(h) := \mathbb{Z}'(h)\mathbf{u} \quad\text{and}\quad \mathbb{X}_\mathbf{u}(h) := (\mathbf{1}_k \vdots \mathbb{Z}'(h)\mathbf{\Gamma_u}),$$

where $\mathbb{Z}(h) := (\mathbf{Z}_{i_1}, \ldots, \mathbf{Z}_{i_k})$ and $\mathbf{1}_k = (1, \ldots, 1)' \in \mathbb{R}^k$, Koenker's result, in the present context, states that $(a_{\boldsymbol{\tau}}^{(n)}, \mathbf{b}_{\boldsymbol{\tau}}^{(n)\prime})' = (\mathbb{X}_{\mathbf{u}}(h))^{-1} \mathbb{Y}_{\mathbf{u}}(h)$ (we just pointed out that, under such conditions, there always exists a quantile hyperplane fitting exactly $k$ observations) is a solution of (2.5) if and only if

$$(3.11) \qquad\qquad -\tau \mathbf{1}_k \leq \boldsymbol{\xi}_{\boldsymbol{\tau}}(h) \leq (1-\tau)\mathbf{1}_k,$$

where

$$(3.12) \qquad \boldsymbol{\xi}_{\boldsymbol{\tau}}(h) := (\mathbb{X}'_{\mathbf{u}}(h))^{-1} \sum_{i \notin h} (\tau - \mathbb{I}_{[r_i < 0]}) \begin{pmatrix} 1 \\ \boldsymbol{\Gamma}'_{\mathbf{u}} \mathbf{Z}_i \end{pmatrix}$$

with $r_i := \mathbf{u}'\mathbf{Z}_i - \mathbf{b}_{\boldsymbol{\tau}}^{(n)} \boldsymbol{\Gamma}'_{\mathbf{u}} \mathbf{Z}_i - a_{\boldsymbol{\tau}}^{(n)}$. Again, this solution is unique if and only if the inequalities in (3.11) are strict; see [16]. As for the constrained case, it follows from the linear programming theory that $(a_{\boldsymbol{\tau}}^{(n)}, \mathbf{c}_{\boldsymbol{\tau}}^{(n)\prime})'$ are the coefficients of a $\boldsymbol{\tau}$-quantile hyperplane if and only if (3.11) holds with $r_i := \mathbf{c}_{\boldsymbol{\tau}}^{(n)\prime} \mathbf{Z}_i - a_{\boldsymbol{\tau}}^{(n)}$ in (3.12) (still with a unique solution when the inequalities are strict).

We stress that no conditions (in particular, no moment conditions) are required here; only, the data points are assumed to be in general position.

3.2. *Equivariance properties.* For the sake of simplicity, results for population quantiles here are stated under Assumption (A); more general statements could be derived, however, by taking into account the possible nonunicity of the resulting $\boldsymbol{\tau}$-quantiles (see Proposition 2.1). It is then easy to check that, with obvious notation, the affine-equivariance property

$$(3.13) \qquad \pi_{\boldsymbol{\tau}\mathbf{Mu}/\|\mathbf{Mu}\|}(\mathbf{MZ} + \mathbf{d}) = \mathbf{M}\pi_{\boldsymbol{\tau}\mathbf{u}}(\mathbf{Z}) + \mathbf{d}$$

holds for any invertible $k \times k$ matrix $\mathbf{M}$ and any vector $\mathbf{d} \in \mathbb{R}^k$. Since, moreover, $\|\tau \mathbf{Mu}\|/\|\mathbf{Mu}\| = \|\tau \mathbf{u}\|$, (3.13) is also compatible with the general equivariance property advocated by [34] in his Definition 2.1. In particular, for translations, we have $\pi_{\boldsymbol{\tau}\mathbf{u}}(\mathbf{Z} + \mathbf{d}) = \pi_{\boldsymbol{\tau}\mathbf{u}}(\mathbf{Z}) + \mathbf{d}$ for any $k$-vector $\mathbf{d}$, which confirms that our concept of multivariate quantiles is *not localized* at any point of the $k$-dimensional Euclidean space; this was not so clear in Section 2 since the center of the unit sphere $\mathcal{S}^{k-1}$ (the origin of $\mathbb{R}^k$) seems to play an important role in their definitions. This is in sharp contrast with other directional quantile contours that are defined with respect to some location center, such as those of [20] (under the terminology *quantile biplots*) and [37].

Note that for any $\tau \in (0, 1)$ and any $\mathbf{u} \in \mathcal{S}^{k-1}$,

$$(3.14) \qquad\qquad \pi_{(1-\tau)\mathbf{u}}(\mathbf{Z}) = \pi_{\tau(-\mathbf{u})}(\mathbf{Z})$$

with the corresponding upper and lower halfspaces exchanged: int $H_{(1-\tau)\mathbf{u}}^{\pm}(\mathbf{Z}) =$ int $H_{\tau(-\mathbf{u})}^{\mp}(\mathbf{Z})$. Clearly, there is no general link between $\pi_{\tau(-\mathbf{u})}(\mathbf{Z})$ and $\pi_{\tau\mathbf{u}}(\mathbf{Z})$ unless the distribution of $\mathbf{Z}$ is centrally symmetric with respect to some point $\boldsymbol{\theta} \in \mathbb{R}^k$.

3.3. *Asymptotic results.* This section derives, under Assumption $(A_n)$ above, strong consistency, asymptotic normality and Bahadur-type representation results for sample $\boldsymbol{\tau}$-quantiles and related quantities.

Under Assumption (A), the population $\boldsymbol{\tau}$-quantiles $(a_{\boldsymbol{\tau}}, \mathbf{b}_{\boldsymbol{\tau}}')'$ and $(a_{\boldsymbol{\tau}}, \mathbf{c}_{\boldsymbol{\tau}}')'$ always are uniquely defined (Proposition 2.1), unlike their sample counterparts $(a_{\boldsymbol{\tau}}^{(n)}, \mathbf{b}_{\boldsymbol{\tau}}^{(n)\prime})'$ and $(a_{\boldsymbol{\tau}}^{(n)}, \mathbf{c}_{\boldsymbol{\tau}}^{(n)\prime})'$; in the sequel, the latter notation will be used for arbitrary sequences of solutions to (2.5) and (2.6), respectively.

Strong consistency of our sample $\boldsymbol{\tau}$-quantiles, namely the fact that $(a_{\boldsymbol{\tau}}^{(n)}, \mathbf{b}_{\boldsymbol{\tau}}^{(n)\prime})'$ converges to $(a_{\boldsymbol{\tau}}, \mathbf{b}_{\boldsymbol{\tau}}')'$ almost surely as $n \to \infty$, holds under Assumption $(A_n)$; this follows, for example, from [13], Section 2.3. Asymptotic normality and Bahadur-type representation results, however, require slightly stronger assumptions. Consider the following reinforcement of Assumption $(A_n)$.

ASSUMPTION $(A_n')$. The observations $\mathbf{Z}_i$, $i = 1, \ldots, n$ are i.i.d. with a common distribution that is absolutely continuous with respect to the Lebesgue measure on $\mathbb{R}^k$, with a density ($f$, say) that has a connected support, admits finite *second*-order moments and, for some constants $C > 0$, $r > k - 2$ and $s > 0$, satisfies

$$(3.15) \quad |f(\mathbf{z}_1) - f(\mathbf{z}_2)| \leq C\|\mathbf{z}_1 - \mathbf{z}_2\|^s \left(1 + \left\|\frac{\mathbf{z}_1 + \mathbf{z}_2}{2}\right\|^2\right)^{-(3+r+s)/2}$$

for all $\mathbf{z}_1, \mathbf{z}_2 \in \mathbb{R}^k$.

Condition (3.15) is very mild. In particular, for $s = 1$, it is satisfied by any continuously differentiable density $f$ for which there exist some constants $C > 0$, $r > k - 2$ and some invertible $k \times k$ matrix $\mathbf{M}$ such that

$$\sup_{\|\mathbf{Mz}\| \geq R} \|\nabla f(\mathbf{z})\| < C(1 + R^2)^{-(r+4)/2}$$

for all $R > 0$. Hence, Assumption $(A_n')$ holds, for example, when the $\mathbf{Z}_i$'s are i.i.d. multinormal or elliptical $t$ with $\nu > 2$ degrees of freedom. Differentiability however is not required, and (3.15) also holds, for instance, for elliptical densities proportional to $\exp(-\|\mathbf{Mz}\|)$ (which are not differentiable at the origin).

As we show in the Appendix (see the proof of Theorem 3.1), Assumption $(A_n')$ implies that the (strictly convex) function $(a, \mathbf{b})' \mapsto \Psi_{\boldsymbol{\tau}}(a, \mathbf{b})$ (see

Definition 2.1) is twice differentiable at $(a_{\boldsymbol{\tau}}, \mathbf{b}_{\boldsymbol{\tau}}')'$, with Hessian matrix

$$\mathbf{H}_{\boldsymbol{\tau}} := \int_{\mathbb{R}^{k-1}} \begin{pmatrix} 1 & \mathbf{x}' \\ \mathbf{x} & \mathbf{x}\mathbf{x}' \end{pmatrix} f((a_{\boldsymbol{\tau}} + \mathbf{b}_{\boldsymbol{\tau}}'\mathbf{x})\mathbf{u} + \boldsymbol{\Gamma}_{\mathbf{u}}\mathbf{x})\, d\mathbf{x}$$

$$= \mathbf{J}_{\mathbf{u}}' \left( \int_{\mathbf{u}^{\perp}} \begin{pmatrix} 1 & \mathbf{z}' \\ \mathbf{z} & \mathbf{z}\mathbf{z}' \end{pmatrix} f((a_{\boldsymbol{\tau}} - \mathbf{c}_{\boldsymbol{\tau}}'\mathbf{z})\mathbf{u} + \mathbf{z})\, d\sigma(\mathbf{z}) \right) \mathbf{J}_{\mathbf{u}} =: \mathbf{J}_{\mathbf{u}}'\mathbf{H}_{\boldsymbol{\tau}}^c \mathbf{J}_{\mathbf{u}},$$

where $\mathbf{u}^{\perp} := \{\mathbf{z} \in \mathbb{R}^k : \mathbf{u}'\mathbf{z} = 0\}$ and $\mathbf{J}_{\mathbf{u}}$ denotes the $(k+1) \times k$ block-diagonal matrix with diagonal blocks 1 and $\boldsymbol{\Gamma}_{\mathbf{u}}$. Strict convexity implies that $\mathbf{H}_{\boldsymbol{\tau}}$ is positive semidefinite. Since, however, for all $\boldsymbol{\tau}$ and $\mathbf{w} := (v_0, \mathbf{v}')' \neq \mathbf{0}$,

$$\mathbf{w}'\mathbf{H}_{\boldsymbol{\tau}}\mathbf{w} = \int_{\mathbb{R}^{k-1}} (v_0 + \mathbf{v}'\mathbf{x})^2 f((a_{\boldsymbol{\tau}} + \mathbf{b}_{\boldsymbol{\tau}}'\mathbf{x})\mathbf{u} + \boldsymbol{\Gamma}_{\mathbf{u}}\mathbf{x})\, d\mathbf{x},$$

$\mathbf{H}_{\boldsymbol{\tau}}$, under Assumption $(A_n')$, is positive definite for all $\boldsymbol{\tau}$.

Letting $\boldsymbol{\xi}_{i,\boldsymbol{\tau}}(a, \mathbf{b}) := -(\tau - \mathbb{I}_{[\mathbf{u}'\mathbf{Z}_i - \mathbf{b}'\boldsymbol{\Gamma}_{\mathbf{u}}'\mathbf{Z}_i - a < 0]})\dot{\mathbf{Z}}_i$ and $\boldsymbol{\xi}_{i,\boldsymbol{\tau}}^c(a, \mathbf{c}) := -(\tau - \mathbb{I}_{[\mathbf{c}'\mathbf{Z}_i - a < 0]})\dot{\mathbf{Z}}_i$, where $\dot{\mathbf{Z}}_i := (1, \mathbf{Z}_i')'$, we have

$$\mathbf{V}_{\boldsymbol{\tau}} := \mathrm{Var}[\mathbf{J}_{\mathbf{u}}'\boldsymbol{\xi}_{1,\boldsymbol{\tau}}(a_{\boldsymbol{\tau}}, \mathbf{b}_{\boldsymbol{\tau}})]$$

$$= \mathbf{J}_{\mathbf{u}}' \begin{pmatrix} \tau(1-\tau) & \tau(1-\tau)\mathrm{E}[\mathbf{Z}'] \\ \tau(1-\tau)\mathrm{E}[\mathbf{Z}] & \mathrm{Var}[(\tau - \mathbb{I}_{[\mathbf{Z}_i \in H_{\boldsymbol{\tau}}]})\mathbf{Z}] \end{pmatrix} \mathbf{J}_{\mathbf{u}}$$

$$= \mathbf{J}_{\mathbf{u}}' \mathrm{Var}[\boldsymbol{\xi}_{1,\boldsymbol{\tau}}^c(a_{\boldsymbol{\tau}}, \mathbf{c}_{\boldsymbol{\tau}})]\mathbf{J}_{\mathbf{u}} =: \mathbf{J}_{\mathbf{u}}'\mathbf{V}_{\boldsymbol{\tau}}^c\mathbf{J}_{\mathbf{u}}.$$

We are then ready to state an asymptotic normality and Bahadur-type representation result for our sample $\boldsymbol{\tau}$-quantile coefficients, which is the main result of this section.

THEOREM 3.1.   *Let Assumption $(A_n')$ hold. Then,*

$$(3.16) \quad \sqrt{n} \begin{pmatrix} a_{\boldsymbol{\tau}}^{(n)} - a_{\boldsymbol{\tau}} \\ \mathbf{b}_{\boldsymbol{\tau}}^{(n)} - \mathbf{b}_{\boldsymbol{\tau}} \end{pmatrix} = -\frac{1}{\sqrt{n}}\mathbf{H}_{\boldsymbol{\tau}}^{-1}\mathbf{J}_{\mathbf{u}}' \sum_{i=1}^n \boldsymbol{\xi}_{i,\boldsymbol{\tau}}(a_{\boldsymbol{\tau}}, \mathbf{b}_{\boldsymbol{\tau}}) + o_{\mathrm{P}}(1)$$

$$(3.17) \qquad\qquad \xrightarrow{\mathcal{L}} \mathcal{N}_k(\mathbf{0}, \mathbf{H}_{\boldsymbol{\tau}}^{-1}\mathbf{V}_{\boldsymbol{\tau}}\mathbf{H}_{\boldsymbol{\tau}}^{-1}) \qquad \text{as } n \to \infty.$$

*Equivalently, writing $\mathbf{P}_k$ for the $(k+1) \times (k+1)$ diagonal matrix with diagonal $(1, -1, \ldots, -1)$,*

$$(3.18) \quad \sqrt{n} \begin{pmatrix} a_{\boldsymbol{\tau}}^{(n)} - a_{\boldsymbol{\tau}} \\ \mathbf{c}_{\boldsymbol{\tau}}^{(n)} - \mathbf{c}_{\boldsymbol{\tau}} \end{pmatrix} = -\frac{1}{\sqrt{n}}\mathbf{P}_k(\mathbf{H}_{\boldsymbol{\tau}}^c)^- \sum_{i=1}^n \boldsymbol{\xi}_{1,\boldsymbol{\tau}}^c(a_{\boldsymbol{\tau}}, \mathbf{c}_{\boldsymbol{\tau}}) + o_{\mathrm{P}}(1)$$

$$(3.19) \qquad\qquad \xrightarrow{\mathcal{L}} \mathcal{N}_{k+1}(\mathbf{0}, \mathbf{P}_k(\mathbf{H}_{\boldsymbol{\tau}}^c)^- \mathbf{V}_{\boldsymbol{\tau}}^c(\mathbf{H}_{\boldsymbol{\tau}}^c)^- \mathbf{P}_k'),$$

*where $(\mathbf{H}_{\boldsymbol{\tau}}^c)^-$ denotes the Moore–Penrose pseudoinverse of $\mathbf{H}_{\boldsymbol{\tau}}^c$. Moreover,*

$$(3.20) \qquad \sqrt{n}(\lambda_{\boldsymbol{\tau}}^{(n)} - \lambda_{\boldsymbol{\tau}}) = \frac{1}{\sqrt{n}} \sum_{i=1}^n (\rho_{\tau}(\mathbf{c}_{\boldsymbol{\tau}}'\mathbf{Z}_i - a_{\boldsymbol{\tau}}) - \lambda_{\boldsymbol{\tau}}) + o_{\mathrm{P}}(1)$$

$$(3.21) \qquad \overset{\mathcal{L}}{\to} \mathcal{N}(0, \mathrm{Var}[\rho_\tau(\mathbf{c}_\tau' \mathbf{Z}_1 - a_\tau)]).$$

As $\rho_\tau(\cdot)$ is a nonnegative function, the distribution of $\sqrt{n}(\lambda_\tau^{(n)} - \lambda_\tau)$ is likely to be skewed for finite $n$ [see (3.20)], which can be partly corrected via a normalizing transformation such as that from [8]. Also, the proof of the above theorem can be easily generalized to derive the asymptotic distribution of vectors of the form $(a_{\tau_1}^{(n)}, \mathbf{b}_{\tau_1}^{(n)\prime}, \ldots, a_{\tau_J}^{(n)}, \mathbf{b}_{\tau_J}^{(n)\prime})'$, $J \in \mathbb{N}_0$.

Theorem 3.1 of course paves the way to inference about $\boldsymbol{\tau}$-quantiles; in particular, it allows to build confidence zones for them. Testing linear restrictions on $\boldsymbol{\tau}$-quantiles coefficients—that is, testing null hypotheses of the form $\mathcal{H}_0 : (a_\tau, \mathbf{b}_\tau')' \in \mathcal{M}(a_0, \mathbf{b}_0, \boldsymbol{\Upsilon}) := \{(a_0, \mathbf{b}_0')' + \boldsymbol{\Upsilon} \mathbf{v} : \mathbf{v} \in \mathbb{R}^\ell\}$ (indexed by some $k$-vector $(a_0, \mathbf{b}_0')'$ and some full-rank $k \times \ell$ matrix $\boldsymbol{\Upsilon}$, $\ell < k$)—can be achieved in the same way as in [25]. Defining and studying such tests requires a detailed investigation of the asymptotic behavior of the constrained estimators

$$(\tilde{a}_\tau^{(n)}, \tilde{\mathbf{b}}_\tau^{(n)\prime})' := \underset{(a, \mathbf{b}')' \in \mathcal{M}(a_0, \mathbf{b}_0, \boldsymbol{\Upsilon})}{\arg\min} \Psi_\tau^{(n)}(a, \mathbf{b}),$$

which is beyond the scope of this work.

## 4. Multivariate quantiles as depth contours.

Turning to the *contour* nature of our multivariate quantiles, we first define the (population and sample) quantile regions and contours that naturally follow from Definitions 2.1 and 2.2 and their empirical counterparts, and state their basic properties. We then establish the strong connections between those regions/contours and the classical Tukey *halfspace depth regions/contours*. Finally, we compare our results with those of Kong and Mizera [20] (Section 4.3) and Wei [37] (Section 4.4).

4.1. *Quantile regions.* The proposed quantile regions are obtained by taking, for some fixed $\tau(=\|\boldsymbol{\tau}\|)$, the "upper envelope" of our $\boldsymbol{\tau}$-quantile hyperplanes. More precisely, for any $\tau \in (0, 1)$, we define our $\tau$-quantile region $R(\tau)$ as

$$(4.1) \qquad R(\tau) := \bigcap_{\mathbf{u} \in \mathcal{S}^{k-1}} \cap \{H_{\tau\mathbf{u}}^+\},$$

where $\cap \{H_{\tau\mathbf{u}}^+\}$ stands for the intersection of the collection $\{H_{\tau\mathbf{u}}^+\}$ of all (closed) upper $(\tau\mathbf{u})$-quantile halfspaces (2.3); for $\tau = 0$, we simply let $R(\tau) := \mathbb{R}^k$. The corresponding $\tau$-quantile contour then is defined as the boundary $\partial R(\tau)$ of $R(\tau)$. At this stage, it is already clear that those $\tau$-quantile regions are *closed* and *convex* (since they are obtained by intersecting closed halfspaces). As we will see below, they are also *nested*: $R(\tau_1) \subseteq R(\tau_2)$ if $\tau_1 \geq \tau_2$.

Empirical quantile regions $R^{(n)}(\tau)$ are obtained by replacing in (4.1) the population quantile halfspaces $H_{\tau\mathbf{u}}^{+}$ with their sample counterparts $H_{\tau\mathbf{u}}^{(n)+}$, yielding, parallel to (4.1),

$$(4.2) \qquad R^{(n)}(\tau) := \bigcap_{\mathbf{u}\in\mathcal{S}^{k-1}} \cap\{H_{\tau\mathbf{u}}^{(n)+}\}$$

for any $\tau \in (0,1)$, with $R^{(n)}(0) := \mathbb{R}^k$. Since they result from intersecting *finitely many* closed halfspaces, these empirical quantile regions are *closed convex polyhedral* sets, the faces of which all are part of some quantile hyperplanes of order $\tau$. Another important property of our empirical regions, which readily follows from the equivariance properties of Section 3.2, is that, for any invertible $k \times k$ matrix $\mathbf{M}$ and any $k$-vector $\mathbf{d}$, using obvious notation,

$$R^{(n)}(\tau; \mathbf{M}\mathbf{Z}_1 + \mathbf{d}, \ldots, \mathbf{M}\mathbf{Z}_n + \mathbf{d}) = \mathbf{M}R^{(n)}(\tau; \mathbf{Z}_1, \ldots, \mathbf{Z}_n) + \mathbf{d}.$$

Similarly, the population regions, in view of (3.13), satisfy the affine-equivariance property $R(\tau; \mathbf{M}\mathbf{Z} + \mathbf{d}) = \mathbf{M}R(\tau; \mathbf{Z}) + \mathbf{d}$ for any such $\mathbf{M}$ and $\mathbf{d}$.

4.2. *Connection with halfspace depth regions.* Recall that the *halfspace* or *Tukey depth* [36] of $\mathbf{z} \in \mathbb{R}^k$ with respect to the probability distribution P is defined as $HD(\mathbf{z}, \mathrm{P}) := \inf\{\mathrm{P}[H] : H$ is a closed halfspace containing $\mathbf{z}\}$. The halfspace depth region $D(\tau)$ of order $\tau \in [0,1]$ associated with P then collects all points of the $k$-dimensional Euclidean space with depth at least $\tau$, that is,

$$(4.3) \qquad D(\tau) = D_{\mathrm{P}}(\tau) := \{\mathbf{z} \in \mathbb{R}^k : HD(\mathbf{z}, \mathrm{P}) \geq \tau\}.$$

Clearly, $D(0) = \mathbb{R}^k$. Also, it is well known (see Proposition 6 in [30], or the proof of Theorem 2.11 in [39] for a more general form) that, for any $\tau > 0$,

$$(4.4) \quad D(\tau) = \cap\{H : H \text{ is a closed halfspace with } \mathrm{P}[\mathbf{Z} \in H] > 1 - \tau\}.$$

The empirical version $D^{(n)}(\tau)$ of $D(\tau)$, as usual, is obtained by replacing, in (4.3) and (4.4), the probability measure P with the empirical measure associated with the observed $n$-tuple $\mathbf{Z}_1, \ldots, \mathbf{Z}_n$ at hand. As shown by the following results, the population halfspace depth regions, under Assumption (A), coincide with the quantile regions $R(\tau)$ defined in (4.1), and so do—almost surely under Assumption $(A_n)$—their empirical counterparts $D^{(n)}(\tau)$, whenever their interior is not empty, with the empirical quantile regions $R^{(n)}(\tau)$ (see the Appendix for the proofs).

THEOREM 4.1. *Under Assumption (A), $R(\tau) = D(\tau)$ for all $\tau \in [0,1)$.*

THEOREM 4.2. *Assume that the $n(\geq k+1)$ data points are in general position. Then, for any $\ell \in \{1, 2, \ldots, n-k\}$ such that $D^{(n)}(\frac{\ell}{n})$ has a nonempty interior, we have that $R^{(n)}(\tau) = D^{(n)}(\frac{\ell}{n})$ for all positive $\tau$ in $[\frac{\ell-1}{n}, \frac{\ell}{n})$.*

Theorem 4.1 of course implies that, under Assumption (A), all results on the halfspace depth regions $D(\tau)$ also apply to the $R(\tau)$ regions. It follows that the $R(\tau)$'s are *compact*; the supremum of all $\tau$'s such that $R(\tau) \neq \varnothing$ belongs to $[1/(k+1), 1/2]$, and takes value $1/2$ if and only if the distribution of $\mathbf{Z}$ is *angularly symmetric*—in the sense that there exists some $k$-vector $\boldsymbol{\theta}$ such that $\frac{\mathbf{Z}-\boldsymbol{\theta}}{\|\mathbf{Z}-\boldsymbol{\theta}\|}$ and $-\frac{\mathbf{Z}-\boldsymbol{\theta}}{\|\mathbf{Z}-\boldsymbol{\theta}\|}$ share the same distribution (see [30] and [32]). This implies that, under Assumption (A), we also may restrict to $\tau \in [0, 1/2]$. As for Theorem 4.2, note that the restriction to halfspace depth regions with nonempty interiors is not really restrictive, since it only applies to flat deepest regions. Another major consequence of this relation between halfspace depth and multivariate quantiles is that our sample multivariate quantiles, just as the traditional univariate ones, completely determine (under the assumptions of Theorem 4.2) the underlying empirical distribution $P_n$—since depth contours do (see [35]). This essentially extends to the population case as well (see, e.g., Section 8 of [20] for a discussion).

Beyond that, Theorems 4.1 and 4.2, by showing that the halfspace depth regions coincide with the upper envelope of *directional* quantile halfspaces, and that the faces of the polyhedral empirical depth contours are parts of empirical quantile hyperplanes, provide depth contours with a straightforward quantile-based interpretation. Above all, Theorem 4.2 brings to the halfspace depth context the extremely efficient computational features of linear programming. This important issue is briefly discussed in Section 5; we refer to [27] for details. See Figure 3 for two- and three-dimensional illustrations.

4.3. *Relation with projection quantiles.* In this section, we discuss the relation of our approach to the results of Kong and Mizera [20] on projection quantiles. These results are somewhat similar to ours, since they also lead to a reconstruction of Tukey's halfspace depth contours. As explained in the Introduction, their $\boldsymbol{\tau}$-quantile is a point in the sample space; denoting by $\tau \mapsto q_\tau^X$ the traditional quantile function associated with the univariate random variable $X$, the $\boldsymbol{\tau}$-quantile $q_{\mathrm{KM};\boldsymbol{\tau}} = q_{\mathrm{KM};\tau\mathbf{u}}$ of a random vector $\mathbf{Z}$ (actually, of its distribution) is defined as $q_\tau^{\mathbf{u}'\mathbf{Z}}\mathbf{u}$, with upper and lower quantile halfspaces

$$(4.5) \qquad H^+_{\mathrm{KM};\tau\mathbf{u}} := \{\mathbf{z} \in \mathbb{R}^k : \mathbf{u}'\mathbf{z} \geq \mathbf{u}'\mathbf{q}_{\mathrm{KM};\tau\mathbf{u}}\}$$

and

$$H^-_{\mathrm{KM};\tau\mathbf{u}} := \{\mathbf{z} \in \mathbb{R}^k : \mathbf{u}'\mathbf{z} < \mathbf{u}'\mathbf{q}_{\mathrm{KM};\tau\mathbf{u}}\},$$

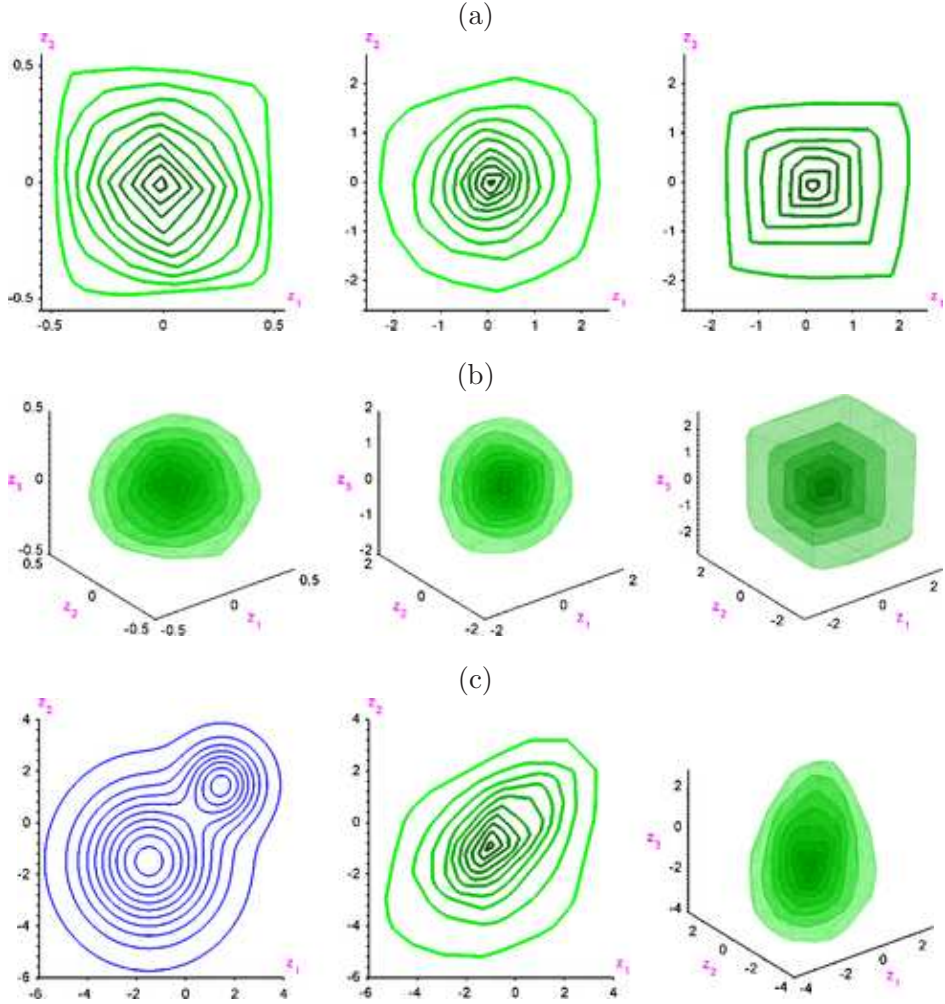

Fig. 3. *Tukey contours $D^{(n)}(\tau)$ (in green) obtained for $n = 449$ from $U([-0.5, 0.5]^k)$, $\mathcal{N}(0,1)^k$, and $t_1^k$ (the products of $k$ independent uniform, standard Gaussian and Cauchy distributions, respectively),* (a) *for $k = 2$ and $\tau \in \{0.01, 0.05, 0.10, 0.15, 0.20, \ldots, 0.45\}$, and* (b) *for $k = 3$ and $\tau \in \{0.05, 0.10, 0.15, 0.20, \ldots, 0.40\}$. For the same $n$ and $\tau$'s, the Tukey depth contours (in green) from the mixtures (with obvious notation) $0.2 \times \mathcal{N}(1.5, 1)^k + 0.8 \times \mathcal{N}(-1.5, 3)^k$ are provided for $k = 2$ and $3$ in* (c)*, along with the density contours (in blue) for $k = 2$. Only the contours falling in the plotting range are displayed.*

respectively, and quantile hyperplane $\pi_{\mathrm{KM};\tau\mathbf{u}} := \{\mathbf{z} \in \mathbb{R}^k : \mathbf{u}'\mathbf{z} = \mathbf{u}'\mathbf{q}_{\mathrm{KM};\tau\mathbf{u}}\}$. Note that those hyperplanes, contrary to ours, are orthogonal to $\mathbf{u}$, so that the relation between $\mathbf{u}$ and $\pi_{\mathrm{KM};\tau\mathbf{u}}$ does not carry any information. Kong

and Mizera show that

$$(4.6) \qquad R_{\mathrm{KM}}(\tau) := \bigcap_{\mathbf{u} \in \mathcal{S}^{k-1}} \{H^+_{\mathrm{KM};\tau\mathbf{u}}\} = D(\tau) \qquad \text{for any } \tau$$

and that

$$(4.7) \qquad R^{(n)}_{\mathrm{KM}}(\tau) = D^{(n)}\left(\frac{\ell}{n}\right) \qquad \text{for any } \tau \in \left[\frac{\ell-1}{n}, \frac{\ell}{n}\right)$$

(see [20] and [26] for different proofs of this latter equality), where $R^{(n)}_{\mathrm{KM}}(\tau)$ stands for the empirical version of $R_{\mathrm{KM}}(\tau)$, obtained by replacing P with the empirical measure $\mathrm{P}_n$ associated with a sample of size $n$.

The results in (4.6) and (4.7) at first sight look pretty equivalent to those of Theorems 4.1 and 4.2, since they also establish a close connection between depth and directional quantiles—here, the Kong and Mizera ones. That connection in (4.6) and (4.7), however, is much less exploitable than in Theorems 4.1 and 4.2. It does provide the faces of the polyhedral empirical depth regions $D^{(n)}(\tau)$ with a neat and interesting quantile interpretation: each face of $D^{(n)}(\tau)$ indeed is part of the Kong and Mizera quantile hyperplane $\pi^{(n)}_{\mathrm{KM};\tau\mathbf{u}_0}$, where $\mathbf{u}_0$ stands for the unit vector orthogonal to that face and pointing to the interior of $D^{(n)}(\tau)$. Unless the depth region $D^{(n)}(\tau)$ is available from some other source, this is not really helpful, though, since, contrary to the collection $\{\pi^{(n)}_{\tau\mathbf{u}}\}$, which is finite, the collection $\{\pi^{(n)}_{\mathrm{KM};\tau\mathbf{u}}\}$, for fixed $\tau$, contains infinitely many hyperplanes (one for each $\mathbf{u} \in \mathcal{S}^{k-1}$). And, since the definition of the upper envelopes $R^{(n)}_{\mathrm{KM}}(\tau)$ of halfspaces $H^{(n)+}_{\mathrm{KM};\tau\mathbf{u}}$ involves an infinite number of such $H^{(n)+}_{\mathrm{KM};\tau\mathbf{u}}$'s, (4.7), contrary to Theorem 4.2, does not readily provide a feasible computation of $D^{(n)}(\tau)$. It is crucial to understand, in that respect, that our quantile halfspaces $H^{(n)+}_{\tau\mathbf{u}}$ are piecewise constant functions of $\mathbf{u}$, in sharp contrast with their Kong and Mizera counterparts $H^{(n)+}_{\mathrm{KM};\tau\mathbf{u}}$: since $\partial H^{(n)+}_{\mathrm{KM};\tau\mathbf{u}}$ is orthogonal to $\mathbf{u}$ for any direction $\mathbf{u}$, there are uncountably many such upper halfspaces in any neighborhood of any fixed direction $\mathbf{u}$, even in the empirical case. To palliate this, Kong and Mizera [20] propose to sample the unit sphere $\mathcal{S}^{k-1}$, which leads to approximate envelopes, that only approximately satisfy (4.7). Moreover, denoting by $\mathbf{U}$ a random vector uniformly distributed over $\mathcal{S}^{k-1}$ (independent of the sample), the probability that the corresponding quantile hyperplane $\pi^{(n)}_{\mathrm{KM};\tau\mathbf{U}}$ contains some face of the Tukey depth contour of order $\tau$ is zero: with probability one, the proposed approximation, thus, fails to recover any of the faces of the actual depth contours. And, for a given sample size, the quality of the approximation deteriorates extremely fast as $k$ increases.

4.4. *Relation to Wei's conditional quantiles.* Another definition of multivariate quantiles, which also extends from location to multiple-output regression, has been proposed by Wei in [37]; see also [38]. Just as Kong and Mizera's projection quantiles and ours, Wei's quantiles are directional quantiles, associated with unit vectors $\mathbf{u} \in \boldsymbol{\mu} + \mathcal{S}^{k-1}$; a center $\boldsymbol{\mu} \in \mathbb{R}^k$ here has to be chosen for the unit sphere—a choice that does have an impact on the final result. Unlike Kong and Mizera's and ours, which are characterized globally, Wei's quantiles, however, are *conditional* ones: the quantiles associated with $\mathbf{u}$ indeed follow from conditional (on $\mathbf{u}$) outlyingness probabilistic characterizations. As a consequence, they are of a local (with respect to $\mathbf{u}$) nature, and their empirical versions therefore unavoidably involve some nonparametric smoothing steps (see Remark 1 on page 399 of [37]). The resulting contours are not convex—hence cannot coincide with depth contours—and strongly depend on the choice of the centering $\boldsymbol{\mu}$.

**5. Computational aspects.** Computational issues in this context are crucial, and we therefore briefly discuss them here. We first restrict to the problem of computing (fixed-$\mathbf{u}$) directional quantiles and related quantities such as the corresponding Lagrange multipliers $\lambda_\tau^{(n)}$ in (3.8c), then consider the computation of (fixed-$\tau$) quantile contours.

5.1. *Computing directional quantiles.* As we have seen in the previous sections, the constrained formulation (2.4) of the definition of our directional quantiles is richer than the unconstrained one (2.1), since it introduces Lagrange multipliers, which bear highly relevant information (that can be exploited for statistical inference; see Section 8). It is therefore natural to focus on the computation of the sample quantiles $(a_\tau^{(n)}, \mathbf{c}_\tau^{(n)\prime})'$ in (2.6) first.

The problem of finding $(a_\tau^{(n)}, \mathbf{c}_\tau^{(n)\prime})'$ can be reformulated as the linear program (P)

$$\min_{(a,\mathbf{c}',\mathbf{r}_+',\mathbf{r}_-')' \in \mathbb{R} \times \mathbb{R}^k \times \mathbb{R}^n \times \mathbb{R}^n} \tau \mathbf{1}_n' \mathbf{r}_+ + (1-\tau) \mathbf{1}_n' \mathbf{r}_-$$

subject to

$$(5.1) \quad \mathbf{u}'\mathbf{c} = 1, \qquad \mathbb{Z}_n'\mathbf{c} - a\mathbf{1}_n - \mathbf{r}_+ + \mathbf{r}_- = \mathbf{0}, \qquad \mathbf{r}_+ \geq \mathbf{0}, \mathbf{r}_- \geq \mathbf{0},$$

where we set $\mathbb{Z}_n := (\mathbf{Z}_1, \ldots, \mathbf{Z}_n)$ and write $\mathbf{r}_+ := (\max(r_1, 0), \ldots, \max(r_n, 0))'$ and $\mathbf{r}_- := (-\min(r_1, 0), \ldots, -\min(r_n, 0))'$. Associated with problem (P) is the dual problem (D)

$$\max_{(\lambda_D, \boldsymbol{\mu}')' \in \mathbb{R} \times \mathbb{R}^n} \lambda_D,$$

subject to

$$\mathbf{1}_n' \boldsymbol{\mu} = 0, \qquad \lambda_D \mathbf{u} + \mathbb{Z}_n \boldsymbol{\mu} = \mathbf{0}_m, \qquad -\tau \mathbf{1}_n \leq \boldsymbol{\mu} \leq (1-\tau)\mathbf{1}_n,$$

where $\lambda_D$ and $\boldsymbol{\mu}$ are the Lagrange multipliers corresponding to the first and second equality constraint in (5.1), respectively. Both (P) and (D) have at least one feasible solution (and therefore also an optimal one). This dual formulation leads to a natural multiple-output generalization of the powerful concept of regression rank scores introduced in [12], allowing for a depth-related form of rank-based inference in this context. This promising line of investigation is not considered here, and left for future research.

We need not worry about the possible nonunicity of the optimal solutions to (P) since, as we have seen in Section 3.3, any sequence of such solutions converges [under Assumption $(A_n)$] to the unique population coefficient vector $(a_{\boldsymbol{\tau}}, \mathbf{c}'_{\boldsymbol{\tau}})'$ almost surely as $n \to \infty$. In practice, one could compute $(a_{\boldsymbol{\tau}}^{(n)}, \mathbf{c}_{\boldsymbol{\tau}}^{(n)'})'$ by means of standard quantile regression of $\mathbf{0}_n$ on $(\mathbf{1}_n | \mathbb{Z}'_n)$ with an extra pseudo-observation consisting of response $C$ and corresponding design row $(0, C\mathbf{u}')$ for some sufficiently large constant $C$, which, in the limit, guarantees that the boundary constraint $\mathbf{u}'\mathbf{c}_{\boldsymbol{\tau}} = 1$ is satisfied; see [2] for another application of the same trick.

Now, since $\lambda_D$ and $\lambda_{\boldsymbol{\tau}}^{(n)}$ are Lagrange multipliers associated with the same constraint, the optimal value $\lambda_D$ of (D) satisfies

$$\lambda_D = n\lambda_{\boldsymbol{\tau}}^{(n)},$$

where, in view of (3.8c), $\lambda_{\boldsymbol{\tau}}^{(n)}$ has a clear meaning. Besides, due to the Strong Duality Theorem, the optimal values of the objective functions in (P) and (D) coincide. Therefore, $\lambda_{\boldsymbol{\tau}}^{(n)}$ is always unique and one has, with $\Psi_{\boldsymbol{\tau}}^{c(n)}$ defined in (2.6),

$$(5.2) \qquad \lambda_{\boldsymbol{\tau}}^{(n)} = \Psi_{\boldsymbol{\tau}}^{c(n)}(a_{\boldsymbol{\tau}}^{(n)}, \mathbf{c}_{\boldsymbol{\tau}}^{(n)}) > 0$$

(except for the rare case of exact fit where $\lambda_{\boldsymbol{\tau}}^{(n)} = 0$), which holds for all optimal solutions to (P) and (D). In other words, $\lambda_{\boldsymbol{\tau}}^{(n)}$ can be obtained from solving (P) as a by-product.

Most importantly, (5.2) allows us to focus on computing our $\boldsymbol{\tau}$-quantiles through the unconstrained problem (2.5) without any loss of generality because we may simply set $\lambda_{\boldsymbol{\tau}}^{(n)} = \Psi_{\boldsymbol{\tau}}^{(n)}(a_{\boldsymbol{\tau}}^{(n)}, \mathbf{b}_{\boldsymbol{\tau}}^{(n)})$. This approach is of course advantageous because it falls directly into the realm of quantile regression, as the problem of finding the sample $\boldsymbol{\tau}$-quantiles in (2.5) can be viewed as looking for standard—that is, single-output—regression quantiles in the regression of $\mathbf{Z_u}$ on the marginals of $\mathbf{Z_u^{\perp}}$ and a constant (in the notation of Section 2).

Needless to say, this interpretation has a large number of implications. Above all, it offers fast, powerful and simple tools for computing sample $\boldsymbol{\tau}$-quantiles (along with the corresponding Lagrange multiplier $\lambda_{\boldsymbol{\tau}}^{(n)}$) in any fixed direction $\mathbf{u}$ and possibly for all $\tau$'s at once, with $\boldsymbol{\tau} = \tau\mathbf{u}$ as usual.

In particular, there is an excellent package for advanced quantile regression analysis in R (see [17]) and the key function for computing quantile regression estimates is also freely available for MATLAB, for example, from Roger Koenker's homepage at http://www.econ.uiuc.edu/~roger/research/rq/rq.html.

5.2. *Computing quantile contours.* As the previous subsection shows that the computation of $H_{\tau\mathbf{u}}^{(n)+}$ is pretty straightforward, we now turn to the problem of aggregating, as efficiently as possible, the information associated with the various fixed-$\tau$ directional quantile halfspaces in order to compute the $R^{(n)}(\tau)$ regions defined in (4.2). The main issue here lies in the proper identification of the finite set of upper quantile halfspaces characterizing $R^{(n)}(\tau)$. This can be achieved efficiently, for any given $\tau \neq \frac{\ell}{n}, \ell \in \{0, 1, \ldots, n\}$, via *parametric linear programming* techniques. By restricting (here and in Figures 1 to 8) to such $\tau$ values, we avoid—without any loss of generality, in view of Theorem 4.2—the problems related with possibly multiple solutions of (P) for integer values of $n\tau$.

For any fixed $\tau \neq \frac{\ell}{n}$, parametric linear programming indeed reveals that, under Assumption $(A_n)$, $\mathbb{R}^k$ almost surely can be segmented into a finite number of nondegenerate cones $\mathcal{C}_i(\tau)$, $i = 1, 2, \ldots, N_C$, such that

$$(a_{\tau\mathbf{u}}^{(n)}, \mathbf{c}_{\tau\mathbf{u}}^{(n)\prime}) = (a_i, \mathbf{c}_i')/\mathbf{t}_i'\mathbf{u}$$

$$\lambda_{\tau\mathbf{u}}^{(n)} = \lambda_i/\mathbf{t}_i'\mathbf{u}$$

$$\mu_{j,\tau\mathbf{u}}^{(n)} = \begin{cases} \mathbf{v}_{ij}'\mathbf{u}/\mathbf{t}_i'\mathbf{u} \in [-\tau, 1-\tau], & \text{if } r_j = 0, \\ -\tau, & \text{if } r_j > 0, \\ 1-\tau, & \text{if } r_j < 0, \end{cases}$$

with $r_j := \mathbf{c}_j'\mathbf{Z}_j - a_j$, for any $\mathbf{u} \in \mathcal{C}_i(\tau) \cap \mathcal{S}^{k-1}$, $i = 1, 2, \ldots, N_C$ and $j = 1, \ldots, n$; see [27] for further details. Each cone $\mathcal{C}_i(\tau)$ then corresponds to one optimal basis $\mathbb{B}_i = \mathbb{B}_{i,\mathbf{u}}$ that uniquely determines constant scalars and vectors $\lambda_i$, $a_i$, $\mathbf{c}_i$, $\mathbf{v}_{ij}$ and $\mathbf{t}_i$ and guarantees that $\mathbf{t}_i'\mathbf{u} > 0$ for any $\mathbf{u} \in \mathcal{C}_i(\tau) \cap \mathcal{S}^{k-1}$. Consequently, each cone $\mathcal{C}_i(\tau)$ corresponds to exactly one quantile hyperplane, and any statistic $S_{\mathbf{u}}$ of the form

$$S_{\mathbf{u}} = g_1(\lambda_{\mathbf{u}}, a_{\mathbf{u}}, \mathbf{c}_{\mathbf{u}})/g_2(\lambda_{\mathbf{u}}, a_{\mathbf{u}}, \mathbf{c}_{\mathbf{u}})$$

is piecewise constant on the unit sphere whenever $g_1(\lambda, a, \mathbf{c})$ and $g_2(\lambda, a, \mathbf{c})$ are homogenous functions of the same order. Figure 4 provides such cones for a bivariate dataset.

It remains to note that we may investigate all the cones $\mathcal{C}_i(\tau)$'s by passing through them counter-clockwise when $k = 2$. In general, we can use the breadth-first search algorithm and always consider all such $\mathcal{C}_i(\tau)$'s that are adjacent to a cone treated in the previous step and have not been considered yet. If $\mathcal{C}_j(\tau)$ and $\mathcal{C}_i(\tau)$ are adjacent cones with point $\mathbf{u}_f$ inside their common face, then $\mathbb{B}_{j,\mathbf{u}_f}$ (and consequently also $\mathbb{B}_{j,\mathbf{u}}$) may be found from the primal feasible basis $\mathbb{B}_{i,\mathbf{u}_f}$ by only a few iterations of the primal simplex algorithm at most.

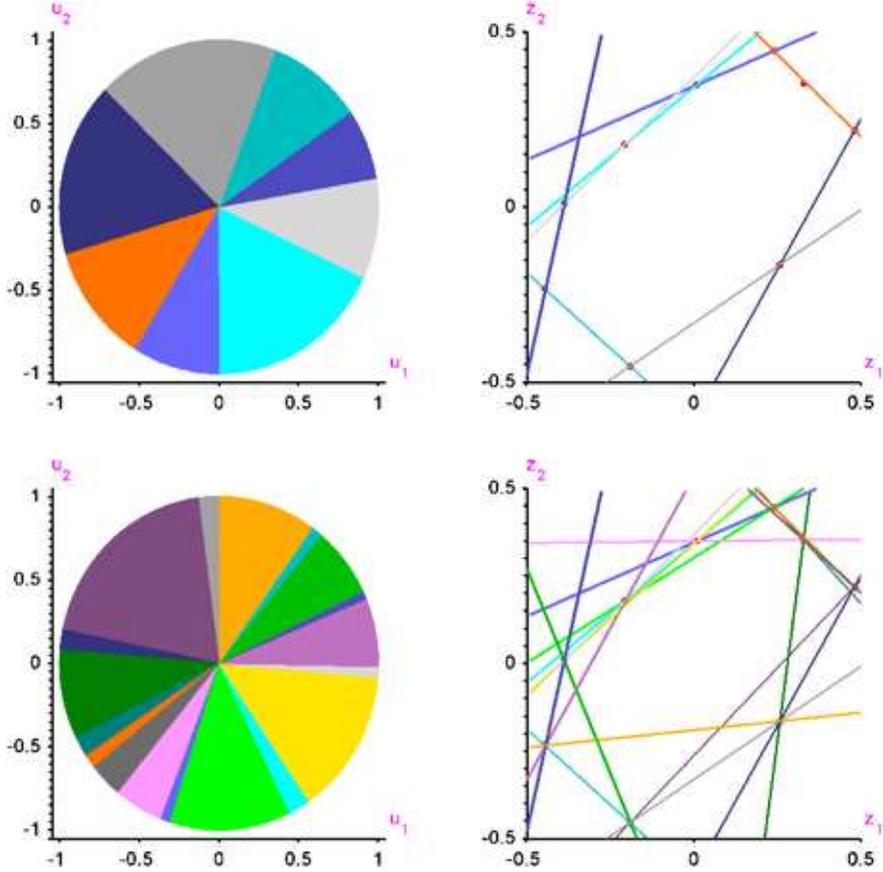

FIG. 4. *The eight cones $\mathcal{C}_i(0.1)$ obtained for $\tau = 0.1$ (top left) and 18 cones $\mathcal{C}_i(0.2)$ obtained for $\tau = 0.2$ (bottom left) via parametric linear programming from the same $n = 9$ points as in Figure 1, along with the corresponding (color matching) $(\tau = 0.1)$-quantile hyperplanes (top right) and $(\tau = 0.2)$-quantile hyperplanes (bottom right).*

Moreover, a careful reading of the proof of Theorem 4.2 reveals (see the remark right after the proof) that a single fixed-$\tau$ collection of quantile hyperplanes $\{\pi_{\tau\mathbf{u}}^{(n)} : \mathbf{u} \in \mathcal{S}^{k-1}\}$ typically contains all hyperplanes relevant for the computation of $k$ consecutive Tukey depth contours. Technical details are provided in [27]. A Matlab implementation of the procedure, which was used to generate all the illustrations in this paper, is available from the authors.

**6. Multiple-output quantile regression.** Our approach to multivariate quantiles also allows to define *multiple-output regression quantiles* enjoying all nice properties of their classical single-output counterparts.

Consider the multiple-output regression problem in which the $m$-variate response $\mathbf{Y} := (Y_1, \ldots, Y_m)'$ is to be regressed on the vector of regressors $\mathbf{X} := (X_1, \ldots, X_p)'$, where $X_1 = 1$ a.s. and the other $X_j$'s are random. In the sequel, we let $\mathbf{X} =: (1, \mathbf{W}')'$, so that $\{(\mathbf{w}', \mathbf{y}') : \mathbf{w} \in \mathbb{R}^{p-1}, \mathbf{y} \in \mathbb{R}^m\} = \mathbb{R}^{p-1} \times \mathbb{R}^m$ is the natural space for considering fitted regression "objects." Multiple-output regression quantiles, in that context, can be obtained by applying Definition 2.1 to the $k$-dimensional random vector $\mathbf{Z} := (\mathbf{W}', \mathbf{Y}')'$, $k = p + m - 1$, *with the important restriction that the direction* $\mathbf{u}$ *should be taken in the response space only*, that is, $\mathbf{u} \in \mathcal{S}_{p-1}^{m-1} := \{\mathbf{0}_{p-1}\} \times \mathcal{S}^{m-1} \subset \mathcal{S}^{k-1}$. This directly yields the following definition.

DEFINITION 6.1.   For any $\boldsymbol{\tau} = \tau \mathbf{u}$, with $\tau \in (0, 1)$ and $\mathbf{u} = (\mathbf{0}_{p-1}', \mathbf{u}_{\mathbf{y}}')' \in \mathcal{S}_{p-1}^{m-1}$, the regression $\boldsymbol{\tau}$-quantile of $\mathbf{Y}$ with respect to $\mathbf{X} = (1, \mathbf{W}')'$ is defined as any element of the collection $\Pi_{\boldsymbol{\tau}}$ of hyperplanes $\pi_{\boldsymbol{\tau}} := \{(\mathbf{w}', \mathbf{y}')' \in \mathbb{R}^{p+m-1} : \mathbf{u}_{\mathbf{y}}'\mathbf{y} = \mathbf{b}_{\boldsymbol{\tau}}'\boldsymbol{\Gamma}_{\mathbf{u}}'(\mathbf{w}', \mathbf{y}')' + a_{\boldsymbol{\tau}}\}$ such that

$$(6.1) \qquad (a_{\boldsymbol{\tau}}, \mathbf{b}_{\boldsymbol{\tau}}')' \in \operatorname*{arg\,min}_{(a, \mathbf{b}')' \in \mathbb{R}^{p+m-1}} \Psi_{\boldsymbol{\tau}}(a, \mathbf{b}),$$

where, denoting by $\boldsymbol{\Gamma}_{\mathbf{u}}$ an arbitrary $(p+m-1) \times (p+m-2)$ matrix such that $(\mathbf{u} \vdots \boldsymbol{\Gamma}_{\mathbf{u}})$ is orthogonal, we let $\Psi_{\boldsymbol{\tau}}(a, \mathbf{b}) := \mathrm{E}[\rho_\tau(\mathbf{u}_{\mathbf{y}}'\mathbf{Y} - \mathbf{b}'\boldsymbol{\Gamma}_{\mathbf{u}}'(\mathbf{W}', \mathbf{Y}')' - a)]$.

Although—similarly as in Definition 2.1—the choice of $\boldsymbol{\Gamma}_{\mathbf{u}}$ has no impact on the directional regression quantile $\pi_{\boldsymbol{\tau}}$, it is here natural to take $\boldsymbol{\Gamma}_{\mathbf{u}}$ of the form

$$\boldsymbol{\Gamma}_{\mathbf{u}} = \begin{pmatrix} \mathbf{I}_{p-1} & \mathbf{0} \\ \mathbf{0} & \boldsymbol{\Gamma}_{\mathbf{u}_{\mathbf{y}}} \end{pmatrix},$$

where $\mathbf{I}_{p-1}$ denotes the $(p-1)$-dimensional identity matrix and the $m \times (m-1)$ matrix $\boldsymbol{\Gamma}_{\mathbf{u}_{\mathbf{y}}}$ is such that $(\mathbf{u}_{\mathbf{y}} \vdots \boldsymbol{\Gamma}_{\mathbf{u}_{\mathbf{y}}})$ is orthogonal. The directional regression quantiles in Definition 6.1 then take the form

$$\pi_{\boldsymbol{\tau}} := \{(\mathbf{w}', \mathbf{y}')' \in \mathbb{R}^{p+m-1} : \mathbf{u}_{\mathbf{y}}'\mathbf{y} = \mathbf{b}_{\boldsymbol{\tau}\mathbf{y}}'\boldsymbol{\Gamma}_{\mathbf{u}_{\mathbf{y}}}'\mathbf{y} + \mathbf{b}_{\boldsymbol{\tau}\mathbf{w}}'\mathbf{w} + a_{\boldsymbol{\tau}}\}$$

with $\mathbf{b}_{\boldsymbol{\tau}} = (\mathbf{b}_{\boldsymbol{\tau}\mathbf{w}}', \mathbf{b}_{\boldsymbol{\tau}\mathbf{y}}')'$. Clearly, an equivalent definition of multiple-output regression quantiles can be obtained by extending Definition 2.2 in the same fashion; see [26].

Now, as in the location case, each quantile hyperplane $\pi_{\boldsymbol{\tau}}$ characterizes a lower (open) and an upper (closed) regression quantile halfspace defined as

$$(6.2) \quad H_{\boldsymbol{\tau}}^- := \{(\mathbf{w}', \mathbf{y}')' \in \mathbb{R}^{p+m-1} : \mathbf{u}_{\mathbf{y}}'\mathbf{y} < \mathbf{b}_{\boldsymbol{\tau}\mathbf{y}}'\boldsymbol{\Gamma}_{\mathbf{u}_{\mathbf{y}}}'\mathbf{y} + \mathbf{b}_{\boldsymbol{\tau}\mathbf{w}}'\mathbf{w} + a_{\boldsymbol{\tau}}\}$$

and

$$(6.3) \quad H_{\boldsymbol{\tau}}^{+} := \{(\mathbf{w}', \mathbf{y}')' \in \mathbb{R}^{p+m-1} : \mathbf{u}_{\mathbf{y}}' \mathbf{y} \geq \mathbf{b}_{\boldsymbol{\tau}\mathbf{y}}' \boldsymbol{\Gamma}_{\mathbf{u}_{\mathbf{y}}}' \mathbf{y} + \mathbf{b}_{\boldsymbol{\tau}\mathbf{w}}' \mathbf{w} + a_{\boldsymbol{\tau}}\},$$

respectively. Most importantly, for fixed $\tau(=\|\boldsymbol{\tau}\|) \in (0,1)$, (multiple-output) $\tau$-*quantile regression regions* are obtained by taking the "upper envelope" of our regression $\boldsymbol{\tau}$-quantile hyperplanes. More precisely, for any $\tau \in (0,1)$, we define regression $\tau$-quantile regions $R_{\mathrm{regr}}(\tau)$ as

$$(6.4) \qquad R_{\mathrm{regr}}(\tau) := \bigcap_{\mathbf{u} \in \mathcal{S}_{p-1}^{m-1}} \cap \{H_{\boldsymbol{\tau}\mathbf{u}}^{+}\}$$

[with corresponding regression quantile contours $\partial R_{\mathrm{regr}}(\tau)$], where $H_{\boldsymbol{\tau}\mathbf{u}}^{+}$ denotes the (closed) upper regression $(\tau\mathbf{u})$-quantile halfspace in (6.3). Unlike the location quantile regions ($p = 1$), regression quantile regions ($p > 1$) may be nonnested—an $m$-dimensional form of the familiar *regression quantile crossing phenomenon*.

Finite-sample versions of all regression concepts above are obtained, similarly as in the location case (Section 2), as the natural sample analogs of the corresponding population concepts; see Figures 5 and 6 for an illustration. From a numerical point of view, Section 5.2, with obvious minor changes, still describes how to compute the resulting regression quantile regions $R_{\mathrm{regr}}^{(n)}(\tau)$, with $m$ and $\mathbf{u}_y$ substituded for $k$ and $\mathbf{u}$, respectively.

The Kong and Mizera projection approach also readily generalizes to the multiple-output regression setting. This issue is briefly addressed in Section 11.3 of [20]; see [26] for a detailed comparison with our approach. As for the conditional quantiles of Wei [37], their regression version shares the same local features as their location counterpart.

## 7. A real data application.
In order to illustrate the implementability and data-analytical power of the concepts we are proposing, we now consider a real data example. Since a thorough case study is beyond the scope of this paper, we only present some very partial results of an investigation of the body girth measurement dataset considered in [14]. That dataset consists of joint measurements of nine skeletal and twelve body girth dimensions, along with weight, height and age, in a group of 247 young men and 260 young women, all physically active. We refer to [14] for details; note, however, that these $n = 507$ observations cannot be considered a random sample representative from any well-defined population, so that the regression quantile contours we are providing below should be taken from a descriptive/illustrative point of view only.

For each gender, taking as regressors a constant term and (with notation $W$) either weight, age, height or the body mass index (defined as BMI:=weight/height$^2$), we considered all $\binom{9+12}{2} = 210$ possible bivariate output regression models, and computed the regression tubes for $\tau = 0.01$, 0.03, 0.10, 0.25 and 0.40, respectively. Three-dimensional pictures of those tubes are not easy to read, and we rather plot, for each of them, a series of five cuts. These cuts were obtained as the intersections of the regression tube under study with hyperplanes of the form $w = w_{(p)}$, where $w_{(p)}$ stands for the (empirical) $p$th quantile of the covariate $W$, $p = 0.10$ (black), 0.30 (blue), 0.50 (green), 0.70 (cyan) and 0.90 (yellow). The results are presented, for women, $Y_1$ the calf maximal girth, and $Y_2$ the thigh maximal girth, with $W$ the weight, age, BMI and height, respectively, in Figure 7.

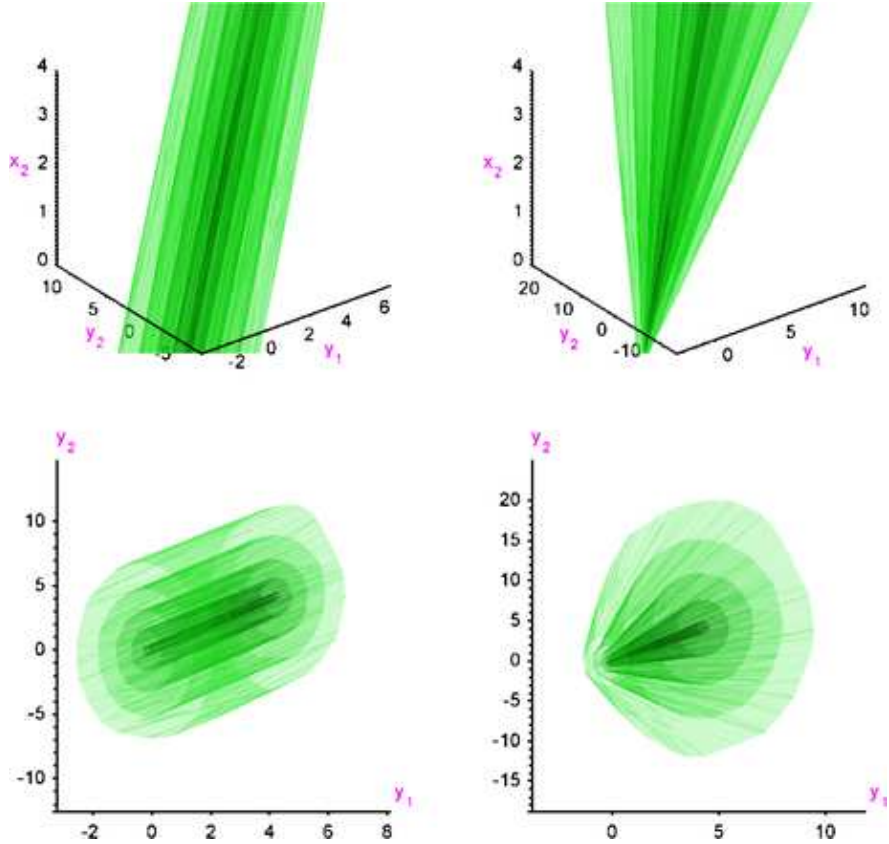

Fig. 5. *Two different views on the regression $\tau$-quantile contours (in green) from $n = 9{,}999$ data points for $\tau \in \{0.01, 0.05, 0.15, 0.30, 0.45\}$ in a homoscedastic $[(Y_1, Y_2)' = (X_2, X_2)' + (\varepsilon_1, \varepsilon_2)';$ left] and a heteroskedastic $[(Y_1, Y_2)' = (X_2, X_2)' + \sqrt{X_2}(\varepsilon_1, \varepsilon_2)';$ right] bivariate-output regression setting, respectively, where $X_2 \sim U([0,4])$, and $\varepsilon_1$ and $\varepsilon_2$ are independent centered Gaussian variables with variances 1 and 9, respectively.*

Results look quite different depending on the choice of regressors. Regression with respect to weight shows a clear positive trend in location (all contours), along with an increasing dispersion, and an evolution of "principal directions," yielding higher variability in calf than in thigh girth for lighter weights (horizontal first "principal direction"), while heavier weights tend to exhibit the opposite phenomenon (vertical first "principal direction"). Quite on the contrary, regression with respect to age apparently does not reveal any location trend: the inner contours almost coincide for all age cuts, and "principal directions" (roughly, parallel to the main bisectors) apparently do not change with age. However, the shapes of outer contours vary quite significantly with age, indicating an increasing (with age) simultaneous variability of both calf and thigh girth largest values. While the results for BMI look very similar to those for weight, a new phenomenon appears when height is the regressor, namely a clear regression effect for some contours (the inner ones) but not for the others, so that the asymmetry structure of the conditional distribution strongly depends on height: the conditional distributions seem much more asymmetric for low values of height than for the higher ones.

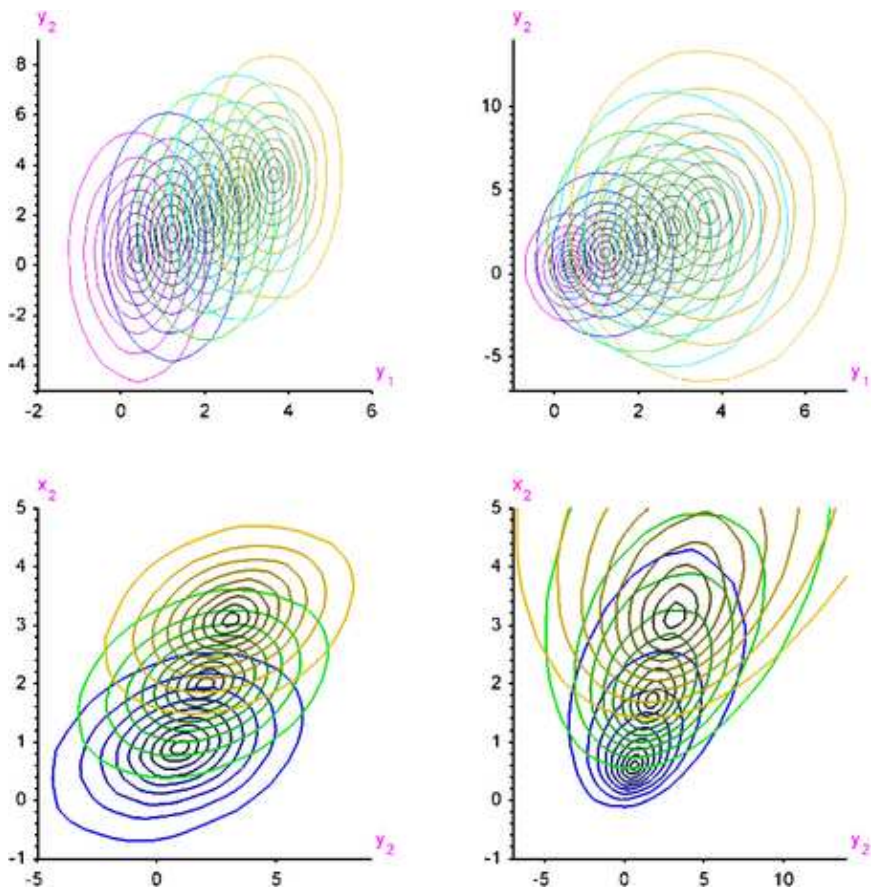

Fig. 6.  *Various cuts of the regression $\tau$-quantile "hypertube" contours from the same two models (left and right, respectively) as in Figure 5 with $n = 9{,}999$ observations. The top plots provide regression $\tau$-quantile cuts, $\tau \in \{0.05, 0.10, 0.15, 0.20, \ldots, 0.45\}$, through 10% (magenta), 30% (blue), 50% (green), 70% (cyan) and 90% (yellow) empirical quantiles of $X_2$; the bottom ones show regression $\tau$-quantile cuts for the same $\tau$ values, and through 25% (blue), 50% (green) and 75% (yellow) empirical quantiles of $Y_1$. Their centers provide information about trend and their shapes and sizes shed light on variability.*

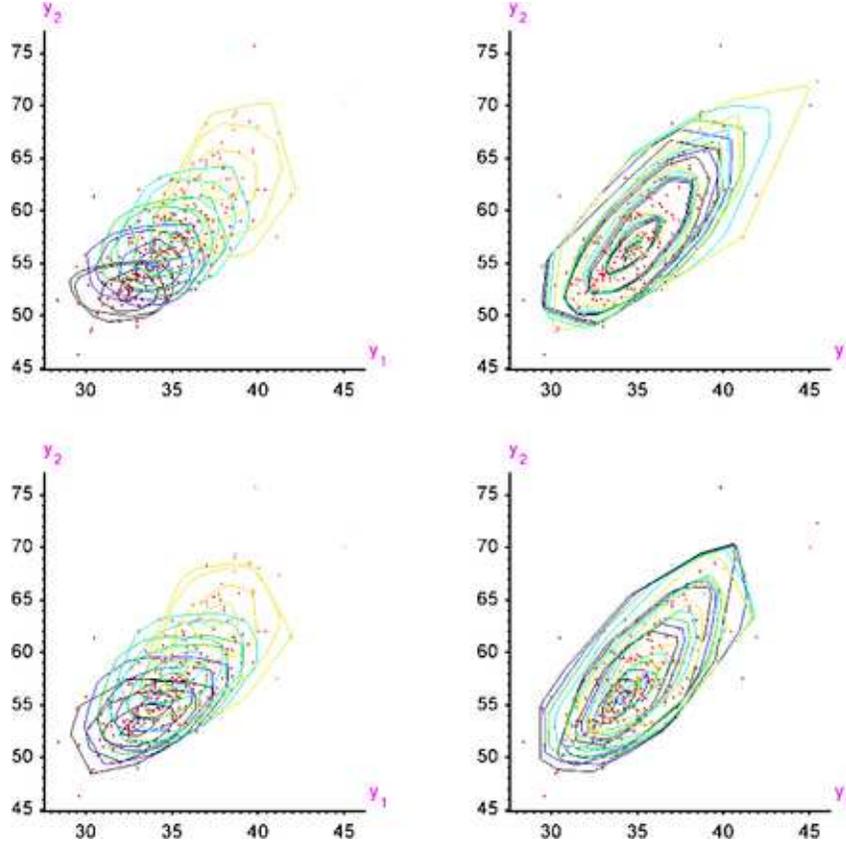

Fig. 7. *Four empirical regression quantile plots from the body girth measurements dataset (women subsample; see [14]). The regression models considered are $(Y_1, Y_2)' = \boldsymbol{\beta}'(1, W)' + (\varepsilon_1, \varepsilon_2)'$, where $Y_1$ is the calf maximum girth and $Y_2$ the thigh maximum girth, while $W$ stands for weight (top left), age (top right), BMI index (bottom left) or height (bottom right). The plots are providing, for $\tau = 0.01$, 0.03, 0.10, 0.25 and 0.40, the cuts of the empirical regression $\tau$-quantile contours, at the empirical $p$-quantiles of the regressors, for $p = 0.10$ (black), 0.30 (blue), 0.50 (green), 0.70 (cyan) and 0.90 (yellow). Data points are shown in red (the lighter the red color, the higher the regressor value).*

These variations in location, scale, shape and asymmetry structures clearly yield a much richer and subtle analysis of the impact of weight/age/BMI/height on those body measurements than any traditional regression method can provide.

**8. Final comments.** This work presents a new concept of multivariate quantile based on $L_1$ optimization ideas and clarifies the quantile nature of halfspace depth contours, while providing an extremely efficient way to compute the latter. The same concept readily allows for an extension of quantile regression to the multiple-output context, thus paving the way to a multiple-output generalization of the many tools and techniques that have been based on the standard (single-output) Koenker and Bassett concept of quantile regression. This final section quickly discusses several open problems, of high practical relevance, that could now be considered.

First of all, Section 6 only very briefly indicates how our multivariate quantiles extend to the context of multiple-output regression; that extension clearly calls for a more detailed study, covering standard asymptotic issues (limiting distributions, Bahadur representations) as well as robustness aspects (breakdown points and influence functions). Nonlinear quantile regression problems also should be addressed via, for instance, local linear methods.

The regression rank score perspectives (associated with linear programming duality) sketched in Section 5.1 also look extremely promising, possibly leading to the development of a full body of multivariate, depth-related methods of rank-based inference.

Finally, as mentioned in the Introduction and in Section 3.1, various concepts introduced in this paper can be quite useful for inference. As an example, note that the symmetry (central, elliptical or spherical) structure of P is reflected in the mappings

$$\mathbf{u} \mapsto \lambda_{\tau\mathbf{u}}\mathbf{u}/\lambda_\tau^{(\infty)} \quad \text{and} \quad \mathbf{u} \mapsto \|\mathbf{c}_{\tau\mathbf{u}}\|\mathbf{u}/c_\tau^{(\infty)},$$

with $\lambda_\tau^{(\infty)} := \sup_{\mathbf{u} \in \mathcal{S}^{k-1}} \lambda_{\tau\mathbf{u}}$ and $c_\tau^{(\infty)} := \sup_{\mathbf{u} \in \mathcal{S}^{k-1}} \|\mathbf{c}_{\tau\mathbf{u}}\|$, as illustrated (with the corresponding empirical quantities, of course) in Figure 8. A test of the hypothesis that the density of $\mathbf{Z}$ is, for example, spherically symmetric thus could be based on [the empirical version $T^{(n)} := T(\mathrm{P}_n)$ of] a functional of the form

$$T(\mathrm{P}) := \int_0^1 \int_{\mathcal{S}^{k-1}} \delta\left(\frac{\lambda_{\tau\mathbf{u}}}{\lambda_\tau^{(\infty)}}, 1\right) d\sigma(\mathbf{u})\, w(\tau)\, d\tau,$$

where $\delta(\cdot, \cdot)$ denotes some distance (such as that of Cramér–von Mises), $w$ some positive weight function over $(0, 1)$, and $\sigma$ the uniform measure over $\mathcal{S}^{k-1}$. Deriving the asymptotic properties of such statistics, however, clearly requires uniform versions of the asymptotic results in Theorem 3.1.

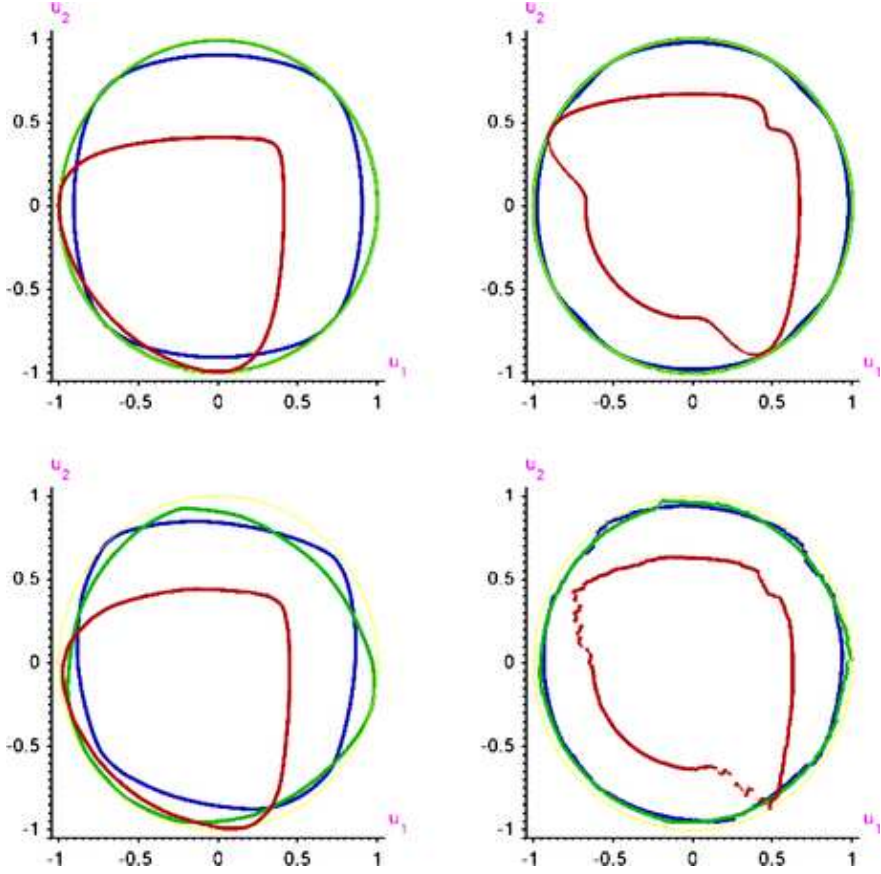

Fig. 8.    *Polar plots of the mappings* $\mathbf{u} \in \mathcal{S}^1 \mapsto \lambda_{\tau\mathbf{u}}^{(n)} \mathbf{u}/(\sup_{\mathbf{v}\in\mathcal{S}^1} \lambda_{\tau\mathbf{v}}^{(n)})$ *(left) and* $\mathbf{u} \in \mathcal{S}^1 \mapsto \|\mathbf{c}_{\tau\mathbf{u}}^{(n)}\| \mathbf{u}/(\sup_{\mathbf{v}\in\mathcal{S}^1} \|\mathbf{c}_{\tau\mathbf{v}}^{(n)}\|)$ *(right), for* $\tau = 0.1$, *and* $n = 49{,}999$ *points (top) [resp.,* $n = 299$ *points (bottom)] drawn from* $\mathcal{N}(0,1)^2$ *(the product of two independent standard Gaussian distributions, in green),* $U([-0.5, 0.5]^2)$ *(the centered bivariate uniform distribution over the unit square, in blue), and* $(\mathrm{Exp}(1) - 1)^2$ *(the product of two independent standard exponential distributions, in red) populations, respectively; see Section* 8. *The resulting shapes clearly reflect the axes of symmetry of the underlying distributions.*

# APPENDIX

Proof of Theorem 3.1.    The quantity $\boldsymbol{\eta}_{i,\tau}(a, \mathbf{b}) := \mathbf{J}_\mathbf{u}' \boldsymbol{\xi}_{i,\tau}(a, \mathbf{b})$ is a subgradient for $(a, \mathbf{b}) \mapsto \rho_\tau(\mathbf{Z}_{i\mathbf{u}} - \mathbf{b}' \mathbf{Z}_{i\mathbf{u}}^\perp - a)$ since, for all $(a, \mathbf{b}')', (a_0, \mathbf{b}_0')' \in \mathbb{R}^k$, we have that

$$\rho_\tau(\mathbf{Z}_{i\mathbf{u}} - \mathbf{b}' \mathbf{Z}_{i\mathbf{u}}^\perp - a) - \rho_\tau(\mathbf{Z}_{i\mathbf{u}} - \mathbf{b}_0' \mathbf{Z}_{i\mathbf{u}}^\perp - a_0) - (a - a_0, \mathbf{b}' - \mathbf{b}_0') \boldsymbol{\eta}_{i,\tau}(a_0, \mathbf{b}_0)$$

$$= (\mathbb{I}_{[\mathbf{u}'\mathbf{Z}_i - \mathbf{b}_0'\boldsymbol{\Gamma}_{i\mathbf{u}}'\mathbf{Z}_i - a_0 < 0]} - \mathbb{I}_{[\mathbf{u}'\mathbf{Z}_i - \mathbf{b}'\boldsymbol{\Gamma}_{i\mathbf{u}}'\mathbf{Z}_i - a < 0]})(\mathbf{u}'\mathbf{Z}_i - \mathbf{b}'\boldsymbol{\Gamma}_{i\mathbf{u}}'\mathbf{Z}_i - a) \geq 0,$$

irrespective of the value of $\mathbf{Z}_i$. Hence, interchanging differentiation and expectation (which is justified in a standard way) shows that $(a, \mathbf{b}')' \mapsto \Psi_{\boldsymbol{\tau}}(a, \mathbf{b})$ (see Definition 2.1) satisfies grad $\Psi_{\boldsymbol{\tau}}(a, \mathbf{b}) = \text{grad } \mathrm{E}[\rho_{\boldsymbol{\tau}}(\mathbf{Z}_{i\mathbf{u}} - \mathbf{b}'\mathbf{Z}_{i\mathbf{u}}^{\perp} - a)] = \mathrm{E}[\boldsymbol{\eta}_{i,\boldsymbol{\tau}}(a, \mathbf{b})]$; see (3.2a) and (3.2b). Therefore,

$$\text{grad } \Psi_{\boldsymbol{\tau}}(a_{\boldsymbol{\tau}} + \Delta_a, \mathbf{b}_{\boldsymbol{\tau}} + \boldsymbol{\Delta}_{\mathbf{b}}) - \text{grad } \Psi_{\boldsymbol{\tau}}(a_{\boldsymbol{\tau}}, \mathbf{b}_{\boldsymbol{\tau}}) - \mathbf{H}_{\boldsymbol{\tau}}(\Delta_a, \boldsymbol{\Delta}_{\mathbf{b}}')'$$

$$= \int_{\mathbb{R}^{k-1}} \int_{a_{\boldsymbol{\tau}} + \mathbf{b}_{\boldsymbol{\tau}}'\mathbf{x}}^{(a_{\boldsymbol{\tau}} + \Delta_a) + (\mathbf{b}_{\boldsymbol{\tau}} + \boldsymbol{\Delta}_{\mathbf{b}})'\mathbf{x}} (f(z\mathbf{u} + \boldsymbol{\Gamma}_{\mathbf{u}}\mathbf{x})$$
$$- f((a_{\boldsymbol{\tau}} + \mathbf{b}_{\boldsymbol{\tau}}'\mathbf{x})\mathbf{u} + \boldsymbol{\Gamma}_{\mathbf{u}}\mathbf{x}))(1, \mathbf{x}')' \, dz \, d\mathbf{x}$$

and Assumption $(A_n')$ yields that

$$\|\text{grad } \Psi_{\boldsymbol{\tau}}(a_{\boldsymbol{\tau}} + \Delta_a, \mathbf{b}_{\boldsymbol{\tau}} + \boldsymbol{\Delta}_{\mathbf{b}}) - \text{grad } \Psi_{\boldsymbol{\tau}}(a_{\boldsymbol{\tau}}, \mathbf{b}_{\boldsymbol{\tau}}) - \mathbf{H}_{\boldsymbol{\tau}}(\Delta_a, \boldsymbol{\Delta}_{\mathbf{b}}')'\|$$

$$\leq C \int_{\mathbb{R}^{k-1}} \left| \int_{a_{\boldsymbol{\tau}} + \mathbf{b}_{\boldsymbol{\tau}}'\mathbf{x}}^{(a_{\boldsymbol{\tau}} + \Delta_a) + (\mathbf{b}_{\boldsymbol{\tau}} + \boldsymbol{\Delta}_{\mathbf{b}})'\mathbf{x}} (|z - (a_{\boldsymbol{\tau}} + \mathbf{b}_{\boldsymbol{\tau}}'\mathbf{x})|^s \|(1, \mathbf{x}')'\|) \right.$$
$$/((1 + \|1/2(z + a_{\boldsymbol{\tau}} + \mathbf{b}_{\boldsymbol{\tau}}'\mathbf{x})\mathbf{u}$$
$$\left. + \boldsymbol{\Gamma}_{\mathbf{u}}\mathbf{x}\|^2)^{(3+r+s)/2}) \, dz \right| d\mathbf{x}$$

$$\leq C \int_{\mathbb{R}^{k-1}} |\Delta_a + \boldsymbol{\Delta}_{\mathbf{b}}'\mathbf{x}| \frac{|\Delta_a + \boldsymbol{\Delta}_{\mathbf{b}}'\mathbf{x}|^s \|(1, \mathbf{x}')'\|}{\|(1, \mathbf{x}')'\|^{3+r+s}} \, d\mathbf{x}$$

$$\leq C \|(\Delta_a, \boldsymbol{\Delta}_{\mathbf{b}}')'\|^{1+s} \int_{\mathbb{R}^{k-1}} \|(1, \mathbf{x}')'\|^{-(r+1)} \, d\mathbf{x} = o(\|(\Delta_a, \boldsymbol{\Delta}_{\mathbf{b}}')'\|)$$

as $\|(\Delta_a, \boldsymbol{\Delta}_{\mathbf{b}}')'\| \to 0$. This shows that $(a, \mathbf{b}')' \mapsto \Psi_{\boldsymbol{\tau}}(a, \mathbf{b})$ is twice differentiable at $(a_{\boldsymbol{\tau}}, \mathbf{b}_{\boldsymbol{\tau}}')'$ with Hessian matrix $\mathbf{H}_{\boldsymbol{\tau}}$. Since, moreover, Assumption $(A_n')$ clearly ensures that $\mathrm{E}[\|\boldsymbol{\eta}_{i,\boldsymbol{\tau}}(a, \mathbf{b})\|^2] < \infty$ for all $(a, \mathbf{b}')' \in \mathbb{R}^k$, Theorem 4 of [25] applies, which establishes (3.16). Of course, (3.17) results from (3.16) by the multivariate CLT.

Recall that, under Assumption (A), the unique solution of (2.1) can be written as $(a_{\boldsymbol{\tau}}, \mathbf{b}_{\boldsymbol{\tau}}')' := (a_{\boldsymbol{\tau}}, (-\boldsymbol{\Gamma}_{\mathbf{u}}'\mathbf{c}_{\boldsymbol{\tau}})')'$, where $(a_{\boldsymbol{\tau}}, \mathbf{c}_{\boldsymbol{\tau}}')'$ denotes the unique solution of (2.4). Similarly, any solution $(a_{\boldsymbol{\tau}}^{(n)}, \mathbf{b}_{\boldsymbol{\tau}}^{(n)'})'$ of (2.5) is related to some solution $(a_{\boldsymbol{\tau}}^{(n)}, \mathbf{c}_{\boldsymbol{\tau}}^{(n)'})'$ of (2.6) via the relation $(a_{\boldsymbol{\tau}}^{(n)}, \mathbf{b}_{\boldsymbol{\tau}}^{(n)'})' = (a_{\boldsymbol{\tau}}^{(n)}, (-\boldsymbol{\Gamma}_{\mathbf{u}}'\mathbf{c}_{\boldsymbol{\tau}}^{(n)})')'$. This allows for rewriting (3.16) as

$$(A.1) \quad \sqrt{n}\, \mathbf{P}_k \mathbf{J}_{\mathbf{u}}' \begin{pmatrix} a_{\boldsymbol{\tau}}^{(n)} - a_{\boldsymbol{\tau}} \\ \mathbf{c}_{\boldsymbol{\tau}}^{(n)} - \mathbf{c}_{\boldsymbol{\tau}} \end{pmatrix} = -\frac{1}{\sqrt{n}} \mathbf{H}_{\boldsymbol{\tau}}^{-1} \mathbf{J}_{\mathbf{u}}' \sum_{i=1}^{n} \boldsymbol{\xi}_{i,\boldsymbol{\tau}}^{c}(a_{\boldsymbol{\tau}}, \mathbf{c}_{\boldsymbol{\tau}}) + o_{\mathrm{P}}(1)$$

as $n \to \infty$. By first premultiplying both sides of (A.1) with $\mathbf{P}_k \mathbf{J}_{\mathbf{u}}$, then using $\boldsymbol{\Gamma}_{\mathbf{u}} \boldsymbol{\Gamma}_{\mathbf{u}}' = \mathbf{I}_k - \mathbf{u}\mathbf{u}'$ [which follows from the orthogonality of $(\mathbf{u}:\boldsymbol{\Gamma}_{\mathbf{u}})$] and

$\mathbf{u}'\mathbf{c}_\tau^{(n)} = 1 = \mathbf{u}'\mathbf{c}_\tau$, we obtain

$$\sqrt{n}\begin{pmatrix} a_\tau^{(n)} - a_\tau \\ \mathbf{c}_\tau^{(n)} - \mathbf{c}_\tau \end{pmatrix} = -\frac{1}{\sqrt{n}}\mathbf{P}_k\mathbf{J}_\mathbf{u}\mathbf{H}_\tau^{-1}\mathbf{J}_\mathbf{u}'\sum_{i=1}^n\boldsymbol{\xi}_{i,\tau}^c(a_\tau,\mathbf{c}_\tau) + o_\mathrm{P}(1)$$

as $n \to \infty$. Lemma A.1 below therefore establishes (3.18). Again, the multi-variate CLT then trivially yields (3.19).

Finally, applying Theorem 6 in [25] [more precisely, applying the version (a) of statement (3.8) in that theorem] with $L = \mathbf{I}_k$ and $c = (a_\tau, \mathbf{b}_\tau')'$ yields

$$(A.2) \qquad \begin{aligned} & n\Psi_\tau^{(n)}(a_\tau, \mathbf{b}_\tau) - n\Psi_\tau^{(n)}(a_\tau^{(n)}, \mathbf{b}_\tau^{(n)}) \\ & - \frac{1}{2n}\sum_{i,j=1}^n\boldsymbol{\xi}_{i,\tau}'(a_\tau,\mathbf{b}_\tau)\mathbf{J}_\mathbf{u}\mathbf{H}_\tau^{-1}\mathbf{J}_\mathbf{u}'\boldsymbol{\xi}_{j,\tau}(a_\tau,\mathbf{b}_\tau) = o_\mathrm{P}(1) \end{aligned}$$

as $n \to \infty$. Note that the third term is clearly $O_\mathrm{P}(1)$ as $n \to \infty$. The result then follows by dividing both sides of (A.2) by $\sqrt{n}$, and by using the identities $\lambda_\tau^{(n)} = \Psi_\tau^{(n)}(a_\tau^{(n)}, \mathbf{b}_\tau^{(n)})$ (see the end of Section 5.1) and $\mathbf{u}'\mathbf{z} - \mathbf{b}_\tau'\boldsymbol{\Gamma}_\mathbf{u}'\mathbf{z} - a_\tau = \mathbf{c}_\tau'\mathbf{z} - a_\tau$ for all $\mathbf{z} \in \mathbb{R}^k$. Since (3.7) entails $\lambda_\tau = \mathrm{E}[\rho_\tau(\mathbf{c}_\tau'\mathbf{Z}_i - a_\tau)]$, the CLT yields (3.21). □

In order to complete the proof of Theorem 3.1, it is sufficient to establish the following lemma.

LEMMA A.1. *The matrix* $\mathbf{G}_\tau := \mathbf{J}_\mathbf{u}(\mathbf{J}_\mathbf{u}'\mathbf{H}_\tau^c\mathbf{J}_\mathbf{u})^{-1}\mathbf{J}_\mathbf{u}'$ *is the Moore–Penrose pseudoinverse of* $\mathbf{H}_\tau^c$*, that is,* $\mathbf{G}_\tau$ *is such that* (i) $\mathbf{G}_\tau\mathbf{H}_\tau^c\mathbf{G}_\tau = \mathbf{G}_\tau$, (ii) $\mathbf{H}_\tau^c\mathbf{G}_\tau\mathbf{H}_\tau^c = \mathbf{H}_\tau^c$, (iii) $(\mathbf{G}_\tau\mathbf{H}_\tau^c)' = \mathbf{G}_\tau\mathbf{H}_\tau^c$ *and* (iv) $(\mathbf{H}_\tau^c\mathbf{G}_\tau)' = \mathbf{H}_\tau^c\mathbf{G}_\tau$.

PROOF. (i) This directly follows from trivial computations. (ii) Let $\mathbf{K}_\mathbf{u}$ be the invertible matrix $(\mathbf{J}_\mathbf{u}\vdots\dot{\mathbf{u}})$, where $\dot{\mathbf{u}} := (0, \mathbf{u}')'$. Clearly, $(\mathbf{H}_\tau^c\mathbf{G}_\tau\mathbf{H}_\tau^c - \mathbf{H}_\tau^c)\mathbf{J}_\mathbf{u} = \mathbf{0}$ and the definition of $\mathbf{H}_\tau^c$ implies that $\dot{\mathbf{u}}$ belongs to the null space of $\mathbf{H}_\tau^c$. Hence, $(\mathbf{H}_\tau^c\mathbf{G}_\tau\mathbf{H}_\tau^c - \mathbf{H}_\tau^c)\mathbf{K}_\mathbf{u} = \mathbf{0}$, which establishes the result. (iii) and (iv) Since $\mathbf{J}_\mathbf{u}'\mathbf{J}_\mathbf{u} = \mathbf{I}_k$, $(\mathbf{G}_\tau\mathbf{H}_\tau^c - \mathbf{H}_\tau^c\mathbf{G}_\tau)\mathbf{J}_\mathbf{u} = \mathbf{J}_\mathbf{u} - \mathbf{H}_\tau^c\mathbf{J}_\mathbf{u}(\mathbf{J}_\mathbf{u}'\mathbf{H}_\tau^c\mathbf{J}_\mathbf{u})^{-1} = \mathbf{0}$; the last equality follows, as in the proof of part (ii), by showing that $(\mathbf{J}_\mathbf{u} - \mathbf{H}_\tau^c\mathbf{J}_\mathbf{u}(\mathbf{J}_\mathbf{u}'\mathbf{H}_\tau^c\mathbf{J}_\mathbf{u})^{-1})'\mathbf{K}_\mathbf{u} = \mathbf{0}$. Now, as we also have that $(\mathbf{G}_\tau\mathbf{H}_\tau^c - \mathbf{H}_\tau^c\mathbf{G}_\tau)\dot{\mathbf{u}} = \mathbf{0}$, we conclude that $(\mathbf{G}_\tau\mathbf{H}_\tau^c - \mathbf{H}_\tau^c\mathbf{G}_\tau)\mathbf{K}_\mathbf{u} = \mathbf{0}$, hence that $\mathbf{G}_\tau\mathbf{H}_\tau^c = \mathbf{H}_\tau^c\mathbf{G}_\tau$. This establishes (iii) and (iv) since both $\mathbf{H}_\tau^c$ and $\mathbf{G}_\tau$ are symmetric. □

PROOF OF THEOREM 4.1. Under Assumption (A), it directly follows from (4.4) that, for any $\tau \in (0, 1)$ (note that Theorem 4.1 trivially holds for $\tau = 0$), $D(\tau) = \cap\{H : H$ is a closed halfspace with $P[\mathbf{Z} \in H] \geq 1 - \tau\}$.

Consequently, by noting that any $H_{\mathrm{KM};\tau\mathbf{u}}^+$, $\mathbf{u} \in \mathcal{S}^{k-1}$ [see (4.5)] satisfies $\mathrm{P}[\mathbf{Z} \in H_{\mathrm{KM};\tau\mathbf{u}}^+] = 1 - \tau$ under Assumption (A), it follows from (4.6) that

$$D(\tau) \subset \cap\{H : H \text{ is a closed halfspace with } P[\mathbf{Z} \in H] = 1 - \tau\}$$

$$\subset \bigcap_{\mathbf{u} \in \mathcal{S}^{k-1}} \{H_{\mathrm{KM};\tau\mathbf{u}}^+\} = D(\tau),$$

which entails that, still under Assumption (A),

(A.3)    $D(\tau) = \cap\{H : H \text{ is a closed halfspace with } \mathrm{P}[\mathbf{Z} \in H] = 1 - \tau\}.$

Now, since (3.2a) [equivalently, (3.5a)] implies that any closed quantile halfspace $H_{\tau\mathbf{u}}^+$, $\mathbf{u} \in \mathcal{S}^{k-1}$, satisfies $\mathrm{P}[\mathbf{Z} \in H_{\tau\mathbf{u}}^+] = 1 - \tau$, (A.3) yields that $D(\tau) \subset R(\tau)$. To show $D(\tau) \supset R(\tau)$, consider an arbitrary closed halfspace $H$ with $\mathrm{P}[\mathbf{Z} \in H] = 1 - \tau$. Then $H = H_{\tau\mathbf{u}}^+$, with

$$\mathbf{u} := \frac{(1/(1-\tau))\mathrm{E}[\mathbf{Z}\mathbb{I}_{[\mathbf{Z} \in H]}] - (1/\tau)\mathrm{E}[\mathbf{Z}\mathbb{I}_{[\mathbf{Z} \in \mathbb{R}^k \setminus H]}]}{\|(1/(1-\tau))\mathrm{E}[\mathbf{Z}\mathbb{I}_{[\mathbf{Z} \in H]}] - (1/\tau)\mathrm{E}[\mathbf{Z}\mathbb{I}_{[\mathbf{Z} \in \mathbb{R}^k \setminus H]}]\|},$$

so that $R(\tau) \subset D(\tau)$; see (3.6) and (A.3) again.    □

PROOF OF THEOREM 4.2.    We start with some remarks on sample halfspace depth regions. By (4.4), $D^{(n)}(\frac{\ell}{n})$, for any $\ell \in \{1, 2, \ldots, n-k\}$, coincides with the intersection of all closed halfspaces containing at least $n - \ell + 1$ observations. Actually, one can restrict to closed halfspaces containing exactly $n - \ell + 1$ observations (see [9], page 1805). Also, it can be shown (see [11]) that $D^{(n)}(\frac{\ell}{n})$—provided that its interior is not the empty set—is bounded by hyperplanes containing at least $k$ points that span a $(k-1)$-dimensional subspace of $\mathbb{R}^k$.

Now, fix $\ell \in \{1, 2, \ldots, n-k\}$ such that $D^{(n)}(\frac{\ell}{n})$ has indeed a nonempty interior. Consider an arbitrary closed halfspace $H$ containing exactly $n - \ell + 1$ data points, among which exactly $k$ ($\mathbf{Z}_i$, $i \in h = \{i_1, \ldots, i_k\}$, say) sit in $\partial H$—and actually span $\partial H$, since the data points are assumed to be in general position. It follows from the results stated in the previous paragraph that $D^{(n)}(\frac{\ell}{n})$, under the assumptions of Theorem 4.2, coincides with the intersection of all such halfspaces.

Letting $s_\tau(n, k, \ell) := (n - k - \ell + 1)\tau + (\ell - 1)(\tau - 1)$, define then

(A.4)    $$\mathbf{u} = \frac{\mathbf{T}_D - s_\tau(n, k, \ell)\mathbf{T}_{on}}{\|\mathbf{T}_D - s_\tau(n, k, \ell)\mathbf{T}_{on}\|},$$

where

$$\mathbf{T}_D := \tau \sum_{\mathbf{Z}_i \in H \setminus \partial H} \mathbf{Z}_i + (\tau - 1) \sum_{\mathbf{Z}_i \notin H} \mathbf{Z}_i \quad \text{and} \quad \mathbf{T}_{on} := \frac{1}{k} \sum_{\mathbf{Z}_i \in \partial H} \mathbf{Z}_i.$$

Taking $\boldsymbol{\Gamma_u}$ as in Definition 2.1, one of course has $\boldsymbol{\Gamma'_u T}_D = s_\tau(n, k, \ell) \boldsymbol{\Gamma'_u T}_{on}$. Hence, writing $(a_h, \mathbf{b}'_h)'$ for the unique solution of

$$\mathbf{u'Z}_i - \mathbf{b'}\boldsymbol{\Gamma'_u}\mathbf{Z}_i - a = 0, \qquad i \in h,$$

we obtain

$$\sum_{i \in \{1,\ldots,n\} \setminus h} (\tau - \mathbb{I}_{[\mathbf{u'Z}_i - \mathbf{b}'_h \boldsymbol{\Gamma'_u}\mathbf{Z}_i - a_h < 0]}) \begin{pmatrix} 1 \\ \boldsymbol{\Gamma'_u}\mathbf{Z}_i \end{pmatrix}$$

$$= \tau \sum_{\mathbf{Z}_i \in H \setminus \partial H} \begin{pmatrix} 1 \\ \boldsymbol{\Gamma'_u}\mathbf{Z}_i \end{pmatrix} + (\tau - 1) \sum_{\mathbf{Z}_i \notin H} \begin{pmatrix} 1 \\ \boldsymbol{\Gamma'_u}\mathbf{Z}_i \end{pmatrix}$$

$$= \begin{pmatrix} s_\tau(n, k, \ell) \\ \boldsymbol{\Gamma'_u T}_D \end{pmatrix} = s_\tau(n, k, \ell) \begin{pmatrix} 1 \\ \boldsymbol{\Gamma'_u T}_{on} \end{pmatrix}.$$

Since [see (3.10)]

$$\frac{1}{k} \mathbb{X}'_{\mathbf{u}}(h) \mathbf{1}_k = \begin{pmatrix} 1 \\ \boldsymbol{\Gamma'_u T}_{on} \end{pmatrix},$$

this implies that, with the same notation as in the end of Section 3.1, we have

$$\boldsymbol{\xi}_{\tau\mathbf{u}}(h) = \frac{s_\tau(n, k, \ell)}{k} \mathbf{1}_k,$$

hence that the subgradient conditions (3.11) are satisfied for any $\tau \in [\frac{\ell-1}{n}, \frac{\ell+k-1}{n}]$. It follows that, for any such $\tau$, $H$ coincides with the upper quantile halfspace $H_{\tau\mathbf{u}}^{(n)+}$ associated with some $\pi_{\tau\mathbf{u}}^{(n)} \in \Pi_{\tau\mathbf{u}}^{(n)}$, where $\mathbf{u}$ is as in (A.4), so that

(A.5) $$R^{(n)}(\tau) := \bigcap_{\mathbf{u} \in \mathcal{S}^{k-1}} \cap \{H_{\tau\mathbf{u}}^{(n)+}\} \subset D^{(n)}\left(\frac{\ell}{n}\right)$$

for any *positive* $\tau \in [\frac{\ell-1}{n}, \frac{\ell}{n})$; one should indeed avoid the value $\tau = 0$ for which $R^{(n)}(\tau)$ is not defined as the upper envelope of quantile halfspaces.

Now, fix $\tau \in (0, \frac{\ell}{n})$. Then, according to (3.9), all upper sample quantile halfspaces $H_{\tau\mathbf{u}}^{(n)+}$ generating $R^{(n)}(\tau)$ contain $P + Z \geq \lceil n(1-\tau) \rceil = n - \lfloor n\tau \rfloor \geq n - \ell + 1$ observations. Hence, $D^{(n)}(\frac{\ell}{n}) \subset R^{(n)}(\tau)$ for any such $\tau$. This, jointly with (A.5), establishes the result. $\square$

Most interestingly, the proof of Theorem 4.2 actually shows that, for any $\tau \in (0, 1)$, the set $\{\pi_{\tau\mathbf{u}}^{(n)} : \mathbf{u} \in \mathcal{S}^{k-1}, \pi_{\tau\mathbf{u}}^{(n)}$ contains $k$ data points$\}$ coincides with the collection of all hyperplanes passing through $k$ observations and cutting off at most $\lfloor n\tau \rfloor$ and at least $\lceil n\tau \rceil - k$ data points. Consequently, as stated at the end of Section 5.2, the set of $\tau$-quantile hyperplanes in all directions provides enough material to compute not only one, but $\min(k + \eta_{n\tau}, \lfloor n\tau \rfloor + 1)$ Tukey depth contours at a time, where $\eta_x$ is one if $x$ is an integer and zero otherwise.

**Acknowledgments.** This research originates in the unpublished Master thesis of Benoît Laine [21] defended at the Université libre de Bruxelles (advisor Marc Hallin), devoted to a pioneering exploration of the subject. We are grateful to Ivan Mizera for encouraging Benoît Laine, inviting him to Edmonton, and sharing with him his expertise in depth problems. We are indebted to Bob Serfling for a most helpful reading of an early stage of our manuscript, and we greatly benefited from inspiring discussions with Roger Koenker, Ivan Mizera, Steve Portnoy and Peter Rousseeuw. The careful reading and constructive comments of two referees and an Associate Editor led to substantial improvements of the manuscript. To all of them we express our warmest thanks.

## REFERENCES

[1] ABDOUS, B. and THEODORESCU, R. (1992). Note on the spatial quantile of a random vector. *Statist. Probab. Lett.* **13** 333–336. MR1160756

[2] BASSETT, G. W., KOENKER, R. and KORDAS, G. (2004). Pessimistic portfolio allocation and Choquet expected utility. *J. Fin. Econometrics* **2** 477–492.

[3] BRECKLING, J. and CHAMBERS, R. (1988). *M*-quantiles. *Biometrika* **75** 761–771. MR0995118

[4] BRECKLING, J., KOKIC, P. and LÜBKE, O. (2001). A note on multivariate M-quantiles. *Statist. Probab. Lett.* **55** 39–44. MR1860190

[5] CHAKRABORTY, B. (2001). On affine equivariant multivariate quantiles. *Ann. Inst. Statist. Math.* **53** 380–403. MR1841143

[6] CHAKRABORTY, B. (2003). On multivariate quantile regression. *J. Statist. Plann. Inference* **110** 109–132. MR1944636

[7] CHAUDHURI, P. (1996). On a geometric notion of quantiles for multivariate data. *J. Amer. Statist. Assoc.* **91** 862–872. MR1395753

[8] CHEN, W. W. and DEO, R. S. (2004). Power transformations to induce normality and their applications. *J. R. Stat. Soc. Ser. B Stat. Methodol.* **66** 117–130. MR2035762

[9] DONOHO, D. L. and GASKO, M. (1992). Breakdown properties of location estimates based on halfspace depth and projected outlyingness. *Ann. Statist.* **20** 1803–1827. MR1193313

[10] DUDLEY, R. M. and KOLTCHINSKII, V. I. (1992). The spatial quantiles. Unpublished manuscript.

[11] FUKUDA, K. and ROSTA, V. (2005). Data depth and maximum feasible subsystems. In *Graph Theory and Combinatorial Optimization* (D. Avis, A. Hertz and O. Marcotte, eds.) 37–67. Springer, New York. MR2180130

[12] GUTENBRUNNER, C. and JUREČKOVÁ, J. (1992). Regression rank scores and regression quantiles. *Ann. Statist.* **20** 305–330. MR1150346

[13] HABERMAN, S. J. (1989). Concavity and estimation. *Ann. Statist.* **17** 1631–1661. MR1026303

[14] HEINZ, G., PETERSON, L. J., JOHNSON, R. W. and KERK, C. J. (2003). Exploring relationships in body dimensions. *J. Statist. Education* **11**. Available at http://www.amstat.org/publications/jse/v11n2/datasets.heinz.html.

[15] HETTMANSPERGER, T. P., NYBLOM, J. and OJA, H. (1992). On multivariate notions of sign and rank. In *L₁-Statistical Analysis and Related Methods (Neuchâtel, 1992)* 267–278. North-Holland, Amsterdam. MR1214838

[16] KOENKER, R. (2005). *Quantile Regression. Econometric Society Monographs* **38**. Cambridge Univ. Press, Cambridge. MR2268657

[17] KOENKER, R. (2007). Quantile regression in R: A vignette. Available at http://cran.r-project.org.

[18] KOENKER, R. and BASSETT, G. J. (1978). Regression quantiles. *Econometrica* **46** 33–50. MR0474644

[19] KOLTCHINSKII, V. (1997). M-estimation, convexity and quantiles. *Ann. Statist.* **25** 435–477. MR1439309

[20] KONG, L. and MIZERA, I. (2008). Quantile tomography: Using quantiles with multivariate data. Preprint.

[21] LAINE, B. (2001). Depth contours as multivariate quantiles: A directional approach. Unpublished Master dissertation (advisor M. Hallin). Univ. Libre de Bruxelles, Brussels.

[22] LIU, R. Y., PARELIUS, J. M. and SINGH, K. (1999). Multivariate analysis by data depth: Descriptive statistics, graphics and inference. *Ann. Statist.* **27** 783–840. MR1724033

[23] MILLER, K., RAMASWAMI, S., ROUSSEEUW, P., SELLARÈS, T., SOUVAINE, D., STREINU, I. and STRUYF, A. (2003). Efficient computation of location depth contours by methods of computational geometry. *Stat. Comput.* **13** 153–162. MR1963331

[24] MIZERA, I. (2002). On depth and deep points: A calculus. *Ann. Statist.* **30** 1681–1736. MR1969447

[25] NIEMIRO, W. (1992). Asymptotics for $M$-estimators defined by convex minimization. *Ann. Statist.* **20** 1514–1533. MR1186263

[26] PAINDAVEINE, D. and ŠIMAN, M. (2009). On directional multiple-output quantile regression. ECARES Working Paper 2009-011.

[27] PAINDAVEINE, D. and ŠIMAN, M. (2010). Computing multiple-output regression quantile regions. Submitted.

[28] ROUSSEEUW, P. J. and RUTS, I. (1996). Algorithm AS 307: Bivariate location depth. *J. Appl. Stat.* **45** 516–526.

[29] ROUSSEEUW, P. J. and HUBERT, M. (1999). Regression depth. *J. Amer. Statist. Assoc.* **94** 388–433. With discussion and a reply by the authors and Stefan Van Aelst. MR1702314

[30] ROUSSEEUW, P. J. and RUTS, I. (1999). The depth function of a population distribution. *Metrika* **49** 213–244. MR1731769

[31] ROUSSEEUW, P. J. and STRUYF, A. (1998). Computing location depth and regression depth in higher dimensions. *Stat. Comput.* **8** 193–203.

[32] ROUSSEEUW, P. J. and STRUYF, A. (2004). Characterizing angular symmetry and regression symmetry. *J. Statist. Plann. Inference* **122** 161–173. Contemporary data analysis: Theory and methods. MR2057920

[33] SERFLING, R. (2002). Quantile functions for multivariate analysis: Approaches and applications. *Statist. Neerlandica* **56** 214–232. Special issue: Frontier research in theoretical statistics, 2000 (Eindhoven). MR1916321

[34] SERFLING, R. (2010). Equivariance and invariance properties of multivariate quantile and related functions, and the role of standardization. *J. Nonparametr. Statist.*, in press.

[35] STRUYF, A. and ROUSSEEUW, P. J. (2005). Halfspace depth and regression depth characterize the empirical distribution. *J. Multivariate Anal.* **69** 135–153. MR1701410

[36] TUKEY, J. W. (1975). Mathematics and the picturing of data. In *Proceedings of the International Congress of Mathematicians (Vancouver, B.C., 1974)* **2** 523–531. Canad. Math. Congress, Montreal. MR0426989

[37] WEI, Y. (2008). An approach to multivariate covariate-dependent quantile contours with application to bivariate conditional growth charts. *J. Amer. Statist. Assoc.* **103** 397–409.

[38] WEI, Y., PERE, A., KOENKER, R. and HE, X. (2005). Quantile regression methods for reference growth charts. *Stat. Med.* **25** 1369–1382. MR2226792

[39] ZUO, Y. and SERFLING, R. (2000). General notions of statistical depth function. *Ann. Statist.* **28** 461–482. MR1790005

[40] ZUO, Y. and SERFLING, R. (2000). Structural properties and convergence results for contours of sample statistical depth functions. *Ann. Statist.* **28** 483–499. MR1790006

E.C.A.R.E.S.
UNIVERSITÉ LIBRE DE BRUXELLES
50, AVENUE F.D. ROOSEVELT, CP114
B-1050 BRUXELLES
BELGIUM
E-MAIL: mhallin@ulb.ac.be
          dpaindav@ulb.ac.be
          Miroslav.Siman@ulb.ac.be
URL: http://homepages.ulb.ac.be/~dpaindav